\documentclass[3p]{elsarticle}
\usepackage{titlesec}
\titleformat{\subsection}
{\normalfont\normalsize\bfseries}{\thesubsection}{1em}{}

\usepackage{stmaryrd}
\usepackage{float}
\usepackage{amsmath,amsfonts,amssymb}
\usepackage{graphicx}
\usepackage{caption}
\usepackage{subcaption}
\usepackage{color}
\usepackage{url}
\usepackage{amsmath}
\usepackage{amsfonts}
\usepackage{amssymb}
\usepackage{amsthm}
\usepackage{tcolorbox}

\usepackage{algpseudocode, algorithm} 
\usepackage{algorithmicx}
\usepackage{algpseudocode}
\usepackage[algo2e]{algorithm2e}

\usepackage{hyperref}
\usepackage{siunitx}
\hypersetup{pdfstartview=FitH}
\newtheorem{theorem}{Theorem}
\newtheorem{remark}[theorem]{Remark}
\newtheorem{definition}[theorem]{Definition}
\newtheorem{lemma}[theorem]{Lemma}

\newtheorem{proposition}[theorem]{Proposition}
\allowdisplaybreaks

\newcommand{\bff}{\boldsymbol{f}}
\newcommand{\bfg}{\boldsymbol{g}}
\newcommand{\bfn}{\boldsymbol{n}}

\newcommand{\bfu}{\boldsymbol{u}}
\newcommand{\bfv}{\boldsymbol{v}}
\newcommand{\bfw}{\boldsymbol{w}}

\newcommand{\bfH}{\boldsymbol{H}}
\newcommand{\bfL}{\boldsymbol{L}}
\newcommand{\bfW}{\boldsymbol{W}}

\newcommand{\RR}{\mathbb{R}}
\newcommand{\NN}{\mathbb{N}}
\newcommand{\PP}{\mathbb{P}}
\newcommand{\HH}{\mathbb{H}}
\newcommand{\WW}{\mathbb{W}}

\newcommand{\jump}[1]{{{\llbracket} #1 {\rrbracket}}}
\begin{document}

\title{A finite element scheme for an optimal control problem on steady Navier-Stokes-Brinkman equations}
\author[a]{Jorge Aguayo\corref{cor}}
\ead{jaguayo@dim.uchile.cl}
\address[a]{Center for Mathematical Modeling, Universidad de Chile, Santiago, Chile}
\author[b]{Julie Merten}
\ead{j.y.merten@rug.nl}
\address[b]{Bernoulli Institute, University of Groningen, Groningen, The Netherlands}
\cortext[cor]{Corresponding author}
\date{\today}

\begin{abstract} 
This paper presents a rigorous finite element framework for solving an optimal control problem governed by the steady Navier-Stokes-Brinkman equations, focusing on identifying a scalar permeability parameter $\gamma$ from local velocity observations. Three different finite element discretization schemes are proposed, and a priori error estimates are proven under appropriate regularity assumptions for each one. A key contribution of this paper is the development of residual-based a posteriori error for both fully discrete and semi-discrete schemes, guiding adaptive mesh refinement to achieve comparable accuracy with fewer degrees of freedom. The method of manufactured solutions is used for numerical experiments to validate the theoretical findings, to demonstrate optimal convergence rates and the effectivity index is evaluated to measure their reliability. The framework offers insights into flow control mechanisms and paving the way for extensions to time-dependent, stochastic, or multiphysics problems. 
\end{abstract}

\maketitle

\pagestyle{myheadings} \thispagestyle{plain}

\section{Motivation}
The optimal control of fluid flows in porous media is a fundamental challenge with far-reaching implications in engineering, environmental science, and biomedicine. From managing subsurface contaminant transport to optimizing drug delivery in vascularized tissues, the ability to precisely control fluid dynamics in complex domains is critical. The Navier-Stokes-Brinkman equations, which generalize the Navier-Stokes equations to account for the effect of porous media with permeability parameter $\gamma$, provide a powerful framework for modeling such systems. However, solving the associated optimal control problems, where the goal is to minimize a cost functional subject to the governing PDE constraints, remains computationally demanding and theoretically intricate.

Although the finite element method (FEM) for the Navier-Stokes equations has been extensively studied, for example \cite{GR86,J16}, the applications of FEM to optimal control problems on elliptic equations have been explored in simpler settings in \cite{RV06,T09}. The authors in \cite{KV09} discussed an optimal control problem for a diffusion-reaction equation with a control on the reaction coefficient, with a continuous Lagrange discretization for states and adjoints and a discontinuous and continuous Lagrange discretization for the control, while \cite{FO22} introduced a posteriori error estimates for the control continuous Lagrange discretization and adding an analysis where the control is implicitly given by a variational discretization, defining a semi-discrete scheme. 

Existing work on Navier-Stokes-Brinkman equations focuses on numerical discretization for the forward problem; see, for example \cite{ARV19,GOS18}, while adaptive FEM for optimal control problems in \cite{AFOQ19,SP24} has been limited to simpler linear elliptic PDEs. In \cite{CMR07}, the authors studied the optimal control problem on Navier-Stokes equations, where the control is a right-hand-side term, obtaining a priori error estimates for a FEM scheme of this problem. 
However, the optimal control of Navier-Stokes-Brinkman systems remains underexplored.


Although the Navier-Stokes-Brinkman equations model flow through porous media through an extension of the Stokes-Brinkman equations \cite{A18,L10}, the articles \cite{A99,AC22} have shown that these equations can model immersed obstacles or domain deformations when the permeability term is sufficiently large in such regions. These models have been applied in the modeling of cardiac valves \cite{P23} and ventricular assist devices \cite{F22}. Our work is based on \cite{ABO23,AO22}, which focuses on a distributed parameter identification for Navier-Stokes equations applied to obstacle detection incorporating the Brinkman's law to model porous media effects, for the purpose of studying the numerical solutions of this problem. 

\begin{tcolorbox}[boxrule=0pt, sharp corners]
    The main novelty of this work consists of the following: 
    \begin{itemize} 
        \item A rigorous adaptive finite element framework analysis for the distributed control of steady Navier-Stokes-Brinkman flows establishing three different discrete schemes: the fully discrete schemes where the control is discretized in a discontinuous Galerkin space with degree $0$ as well as in continuous Lagrange space of degree $1$, and the variational scheme, also called semi-discrete scheme.
        \item A priori and a posteriori residual error estimates, based on the estimators defined in \cite{OWA94,V96}, to ensure the reliability and global efficiency of our a posteriori estimator for our three schemes, and the local efficiency of the estimator for the semi-discrete scheme. 
        \item Numerical test to confirm that the a priori error bounds are reached when the optimal control problem is defined for global or local observations of the fluid velocity. The effectivity index $\theta=$(error estimate/true error) which ideally converges to $1$ as the mesh size becomes smaller proves our error estimators to be asymptotically exact\cite{BR79} in most cases.
        \item An adaptive refinement strategy capable of improving the numerical results obtained in \cite{AO22} for the detection of obstacles immersed in fluids.
    \end{itemize}
\end{tcolorbox}

Numerical experiments are conducted using \texttt{FEniCS} \cite{fenics}, a flexible and efficient platform widely used in the scientific community, to validate the proposed framework, demonstrating its effectiveness and practicality. By leveraging \texttt{FEniCS}, we can focus on building a sophisticated algorithm while minimizing the overhead of low-level implementation details. 

The article is structured as follows: Section \ref{S2} introduces the model problem and its optimal control formulation. Section \ref{S3} discusses the finite element discretization scheme for the model problem and the optimal control problem, where cases for discontinuous Galerkin, continuous Lagrange and semi-discrete scheme are addressed. Section \ref{S4} introduces the a posteriori estimator for the adaptive mesh refinement strategy and the effectivity index. Section \ref{S5} shows the computation results of various numerical experiments. And finally, Section \ref{S6} summarizes the core findings and novel contributions, with an outlook for future studies. 

\section{The model problem}\label{S2}

Consider $n\in\{2,3 \}  $ and a non-empty bounded domain $\Omega\subseteq\RR^{n}$. The Lebesgue measure of $\Omega$ is denoted by $\left\vert \Omega\right\vert $, which extends to spaces of lesser dimension. The norm and semi-norms for the Sobolev spaces $W^{m,p}(\Omega)$ are denoted by $\left\Vert \cdot\right\Vert _{m,p,\Omega}$ and $\left\vert \cdot\right\vert _{m,p,\Omega}$, respectively. For $p=2$, the norm, semi-norms, inner product and duality pairing of the space $W^{m,2}(\Omega)=H^{m}(\Omega)$ are denoted by $\left\Vert \cdot\right\Vert _{m,\Omega}$, $\left\vert \cdot\right\vert _{m,\Omega}$, $(\cdot,\cdot)_{m,\Omega}$ and $\left\langle \cdot,\cdot\right\rangle _{m,\Omega}$, respectively. For $m=0$, the inner product is denoted by $(\cdot,\cdot)_\Omega$. In addition, $\mathcal{C}^{m}(\Omega)$ and $\mathcal{C}^{\infty}(\Omega)$ denote the space of functions with continuous derivatives $m$ and all continuous derivatives, respectively. For $\Omega_{1}$ and $\Omega_{2}$ two open subsets of $\RR^{n}$, we denote $\Omega_{1}\Subset\Omega_{2}$ when there exists a compact set $K$ such that $\Omega_{1}\subseteq K\subseteq\Omega_{2}$.

The spaces $\HH^{m}(\Omega)$, $\WW^{m,p}(\Omega)$, $\bfH^{m}(\Omega)$, $\bfw^{m,p}(\Omega)$, $\boldsymbol{C}^{m}(\Omega)$, and $\boldsymbol{C}^{\infty}(\Omega)$ are defined by $\HH^{m}(\Omega)=[  H^{m}(\Omega)]  ^{n\times n}$, $\WW^{m,p}(\Omega)=[  W^{m,p}(\Omega)]  ^{n\times n}$, $\bfH^{m}(\Omega)=[  H^{m}(\Omega)]  ^{n}$, $\bfW^{m,p}(\Omega)=[  W^{m,p}(\Omega)]  ^{n}$, $\boldsymbol{C}^{m}(\Omega)=[  \mathcal{C}^{m}(\Omega)]  ^{n}$ and $\boldsymbol{C}^{\infty}(\Omega)=[  \mathcal{C}^{\infty}(\Omega)] ^{n}$. The notation for norms, semi-norms and inner products will be extended from $W^{m,p}(\Omega)$ or $H^{m}(\Omega)$. Given $\boldsymbol{A},\boldsymbol{B}\in\RR^{n\times n}$, $a_{ij}$ denotes the entry in the $i-$th row and $j- $th column of matrix $\boldsymbol{A}$, $\boldsymbol{A}^{T}$ denotes the transpose matrix of $\boldsymbol{A}$, $\operatorname{tr}(\boldsymbol{A})$ denotes the trace of $\boldsymbol{A}$ and $\boldsymbol{A}:\boldsymbol{B}$ denotes the inner product of $\RR%
^{n\times n}$ given by%
\[
\boldsymbol{A}:\boldsymbol{B}=\operatorname{tr}(\boldsymbol{AB}^{T}%
)=\sum_{i=1}^{n}\sum_{j=1}^{n}a_{ij}b_{ij}\text{,}%
\]
$\operatorname{sym}\boldsymbol{A}=\dfrac{1}{2}(\boldsymbol{A}+\boldsymbol{A}^{T})$ and $\operatorname{skew}\boldsymbol{A}=\dfrac{1}{2}(\boldsymbol{A}-\boldsymbol{A}^{T})$. The identity matrix is denoted by $\boldsymbol{I}\in\RR^{n\times n}$. Analogously, given $\boldsymbol{a},\boldsymbol{b}\in\RR^{n}$, $a_{j}$ denotes the $j-$th component of the vector $\boldsymbol{a}$ and $\boldsymbol{a}\cdot\boldsymbol{b}$ denotes the inner product of $\RR^{n}$.

\subsection{An optimal control problem}
First, we introduce some notation and definitions that allow us to simplify the definition of the optimal control problem to be analyzed.
\begin{definition}
We define the spaces $H=\bfH^{1}(\Omega)=\left\{  \bfv\in\bfH^{1}(\Omega)\text{ } \mid\text{ } \bfv=\boldsymbol{0}\text{ on }\partial\Omega\right\}$ and $Q=L_{0}^{2}(\Omega):=\left\{  q\in L^{2}(\Omega)\text{ } \mid\text{ } (q,1)=0\right\}$. We denote the norm for the space $V:=  V $as $\Vert(v,q)\Vert=\left(\vert \bfv\vert _{1,\Omega}%
^{2}+\Vert q\Vert_{0,\Omega}^{2}\right)  ^{1/2}$
\end{definition}

\begin{definition}
For $\nu>0$, we define for all $\bfu,\bfv,\bfw\in H$ and $q\in Q$.

\begin{enumerate}
\item $a(\bfu,\bfv)=\nu(\nabla\bfu,\nabla\bfv)_\Omega$.
\item $b(\bfv,q)=-(q,\operatorname{div}\bfv)_\Omega$.
\item $c(\bfu,\bfv,\bfw)=((\nabla\bfu)\bfv,\bfv)_\Omega$.
\end{enumerate}
\end{definition}

Consider $\mathcal{A}=\{\gamma\in L^{2}(\Omega)$ $\mid$ $a\leq\gamma\leq b$ in $\Omega\}$, for $a,b\in\RR$ such that $0\leq a<b$. Given a constant $\alpha>0$ and a non-empty open subset $\omega\subset\Omega$, the optimal control problem to be studied is given by
\begin{align}
\text{minimize\quad}   &  J(\gamma)=\dfrac{1}{2}\Vert\bfu-\bfu_{0}\Vert_{0,\omega}^{2}+\dfrac{\alpha}{2}\Vert\gamma -\gamma_{0}\Vert_{0,\Omega}^{2}\label{min}\\
\text{subject to\quad} &  (\forall(\bfw,r)\in V)\quad a(\bfu,\bfw)-b(\bfw,p)+b(\bfu,r)+c(\bfu,\bfu,\bfw)+(\gamma\bfu,\bfv)_\Omega=(\bff,\bfw)_\Omega\label{VF}\\ 
                       &  \gamma\in\mathcal{A},\nonumber
\end{align}
where $\bfu_{0}\in\bfL^{2}(\omega)$, $\gamma_{0}\in L^{2}(\Omega)$ and $\bff\in\bfL^{2}(\Omega)$.

\begin{lemma}\label{bounds-1}
There exist positive constants $\beta$ and $\delta$ such that for all $\bfu,\bfv,\bfw\in H$ and $q\in Q$.

\begin{enumerate}
\item $\vert a(\bfu,\bfv)\vert \leq\nu\vert\bfu\vert_{1,\Omega}\vert\bfv\vert_{1,\Omega}$.
\item $\vert b(\bfv,q) \vert \leq\sqrt{n}\vert\boldsymbol{v}\vert_{1,\Omega}\Vert q\Vert_{0,\Omega}$.
\item $\vert c(\bfu,\bfv,\bfw) \vert\leq\beta\vert\bfu \vert_{1,\Omega}\vert\bfv\vert_{1,\Omega}\vert\bfw\vert_{1,\Omega}$.
\item $\sup\limits_{\bfv\in H\setminus\{\mathbf{0}\}}\dfrac {b(\bfv,q)_\Omega}{\vert \bfv\vert _{1,\Omega}} \geq\delta\Vert q\Vert_{0,\Omega}$.
\end{enumerate}
\end{lemma}

\begin{theorem}\label{f-bound}
Let $\gamma\in\mathcal{A}$. If there exists a constant $C\in(0,1)$ such that $\vert(\bff,\bfv)_{0,\Omega}\vert\leq C\dfrac{\nu^{2}}{\beta}\vert\bfv\vert_{1,\Omega}$ for all $\bfv\in H$, then \eqref{VF} has an unique solution $(\bfu,p)\in  V $ and holds $\vert\bfu\vert_{1,\Omega}\leq C\dfrac{\nu}{\beta}$. Furthermore, there exists $C>0$ such that%
\begin{align*}
\vert \bfu\vert _{1,\Omega}  &  \leq C\Vert\bff\Vert _{0,\Omega},\\
\Vert p\Vert_{0,\Omega}  &  \leq C(\Vert \bff\Vert_{0,\Omega}+\Vert \bff\Vert _{0,\Omega}^{2}).
\end{align*}
\end{theorem}

\begin{proof}
See Chapter 9 in \cite{G11}.
\end{proof}

In what follows, we assume that $\bff$ fulfills the hypothesis of Theorem \ref{f-bound}.
\begin{definition}
We denote the control-to-state map as $A:\mathcal{A}\rightarrow  V $, where $A(\gamma)=(\bfu,p)$ is the solution of \eqref{VF}.
\end{definition}

\begin{definition}
Given $\gamma\in\mathcal{A}$ and $(\bfu,p)=A(\gamma)$, we define the adjoint equation as the variational formulation given by

\qquad\qquad\textit{Find $(\bfv,q)\in V$ such that for all $(\bfw,r)\in V$}
\begin{align}
a(\bfv,\bfw)-b(\boldsymbol{r},q)+b(\bfv,r)+c(\bfu,\bfw,\bfv)+c(\bfw,\bfu,\bfv)+(\gamma\bfv,\bfw)_\Omega=(\bfu-\bfu_{0},\bfw)_{\omega}.\label{AVF}%
\end{align}

\end{definition}

\begin{lemma}\label{adjoint}
Given $\gamma\in\mathcal{A}$, the adjoint problem \eqref{AVF} has an unique solution $(\bfv,q)\in  V $ such that%
\begin{align*}
\vert\bfv\vert_{1,\Omega}  &  \leq C\Vert\bfu-\bfu_{0}\Vert_{0,\omega},\\
&  \leq C(\Vert\bfu_{0}\Vert_{0,\omega}+\Vert\bff\Vert_{0,\Omega})\\
\Vert q\Vert_{0,\Omega}  &  \leq C(\Vert\bfu-\bfu_{0}\Vert_{0,\omega}+\vert\bfv\vert_{1,\Omega}(1+\Vert\bff\Vert_{0,\Omega})\\
&  \leq C\Vert\bff\Vert_{0,\Omega}(\Vert\bfu_{0}\Vert_{0,\omega}+\Vert\bff\Vert_{0,\Omega})
\end{align*}
\end{lemma}
\begin{proof}
It is a direct consequence of the Babuzka-Brezzi theory. We omit the details.
\end{proof}

\begin{definition}
We denote the control-to-adjoint-state map as $B:\mathcal{A}\rightarrow  V $, where $B(\gamma)=(\bfv,q)$ is the solution of \eqref{AVF}.
\end{definition}

\begin{lemma}\label{h2}
Let $\gamma\in\mathcal{A}$, and $(\bfu,p),(\bfv,q)\in  V $ such that $A(\gamma)=(\bfu,p)$ and $B(\gamma)=(\bfv,q)$. Then, $\bfu,\bfv\in\bfH^{2}(\Omega)$ and $p,q\in H^{1}(\Omega)$. Furthermore, the following estimates are verified
\begin{align*}
\vert\bfu\vert_{2,\Omega}+\Vert p\Vert_{1,\Omega}  &  \leq C(\Vert\bff\Vert_{0,\Omega}+\Vert\bff\Vert_{0,\Omega}^{2}),\\
\vert\bfv\vert_{2,\Omega}+\Vert q\Vert_{1,\Omega}  &  \leq C\Vert\bff\Vert_{0,\Omega}(\Vert\bff\Vert_{0,\Omega}+\left\Vert \bfu_{0}\right\Vert _{0,\omega}),
\end{align*}
for a positive constant $C>0$ independent on $\bfu$, $\bfv$, $\bff$ and $\bfu_{0}$.
\end{lemma}

\begin{proof}
See \cite{G92}.
\end{proof}

\begin{theorem}\label{existence}
The optimal control problem \eqref{min} admits at least one global solution $(\boldsymbol{\bar{u}},\bar{p},\bar{\gamma})\in  V \times\mathcal{A}$.
\end{theorem}
\begin{proof}
See Theorem 10 in \cite{ABO23}.
\end{proof}

Now we present some important properties of the cost functional $J$ and the maps $A$ and $B$.

\begin{lemma}\label{lipschitz}
There exist positive constants $C_1$ and $C_2$ such that for all $\gamma_{1},\gamma_{2}\in\mathcal{A}$
\begin{enumerate}
\item $\vert A(\gamma_{1})-A(\gamma_{2})\vert\leq C\Vert\gamma_{1}-\gamma_{2}\Vert_{0,\Omega}$.
\item $\vert B(\gamma_{1})-B(\gamma_{2})\vert\leq C\Vert\gamma_{1}-\gamma_{2}\Vert_{0,\Omega}$.
\end{enumerate}
\end{lemma}

\begin{proof}
See Lemmas 11 and 20 in \cite{ABO23}.
\end{proof}

\begin{proposition}
The functional $J$ and the operators $A$ and $B$ and twice Fr\'{e}chet-differentiable with respect to the $L^{\infty}(\Omega)$ topology. Given $\varphi_{1},\varphi_{2}\in L^{\infty}(\Omega)$ and $\gamma\in\mathcal{A}$, the derivatives of $A$ have the following properties

\begin{enumerate}
\item $A^{\prime}(\gamma)[\varphi_{1}]=(\bfu^{\prime},p^{\prime})$
is the solution of
\[
(\forall(\bfv,q)\in  V )\quad a(\bfu^{\prime},\bfv)-b(\bfv,p^{\prime})+b(\bfu^{\prime},q)+c(\bfu^{\prime},\bfu,\bfv)+c(\bfu,\bfu^{\prime},\bfv)+(\gamma\bfu^{\prime},\bfv)_\Omega=-(\varphi_{1}\bfu,\bfv)_\Omega,
\]
where $A(\gamma)=(\bfu,p)$.

\item $A^{\prime\prime}(\gamma)[\varphi_{1},\varphi_{2}]=(\bfu^{\prime\prime},p^{\prime\prime})$ is the solution of%
\[
(\forall(\bfv,q)\in V )\quad a(\bfu^{\prime\prime},\bfv)-b(\bfv,p^{\prime\prime})+b(\bfu^{\prime\prime},q)+c(\bfu^{\prime\prime},\bfu,\bfv)+c(\bfu,\bfu^{\prime\prime},\bfv)+(\gamma\bfu^{\prime\prime},\bfv)_\Omega=-(\varphi_{1}\bfu_{2}^{\prime}+\varphi_{2}\bfu_{1}^{\prime},\bfv)_\Omega,
\]
where $A^{\prime}(\gamma)[\gamma_{1}]=(\bfu_{1}^{\prime},p_{1}^{\prime})$ and $A^{\prime}(\gamma)[\gamma_{2}]=(\bfu_{2}^{\prime},p_{2}^{\prime})$.
\end{enumerate}
\end{proposition}
\begin{proof}
See Theorems 11 and 18 in \cite{ABO23}.
\end{proof}

\begin{definition}
Consider $\varepsilon>0$. Given $\gamma^\ast\in\mathcal{A}$, we denote $\mathcal{A}_{\varepsilon}(\gamma^\ast)=\{ \varphi \in \mathcal{A} ~ \mid ~ \Vert \varphi - \gamma^\ast \Vert_{0,\Omega}\leq \varepsilon \} $. $\gamma^\ast$ is a local solution of \eqref{min} if and only if there exists $\varepsilon >0$ such that $J(\gamma)\geq J(\gamma^\ast)$ for all $\gamma \in \mathcal{A}_{\varepsilon}(\gamma^\ast)$
\end{definition}

The following results allow us to deduce the first and second-order optimality conditions, together with some qualities that every local solution of \eqref{min} must satisfy. The following is the first-order optimality condition.
\begin{proposition}
If $\bar{\gamma}\in\mathcal{A}$ is a local solution of \eqref{min}, then the following inequality holds
\begin{equation}
(\forall\varphi\in\mathcal{A})\quad J^{\prime}(\bar{\gamma})[\varphi-\bar{\gamma}]\geq0. \label{first-order}
\end{equation}
\end{proposition}

The expression $J^{\prime}(\gamma)[\varphi]$ can be rewritten as
\begin{align*}
J^{\prime}(\gamma)[\varphi]  &  =(\bfu-\bfu_{0},\bfu^{\prime})_{\omega}+\alpha(\gamma-\gamma_{0},\varphi)_\Omega\\
                             &  =(\varphi\bfu,\bfv)_\Omega +\alpha(\gamma-\gamma_{0},\varphi)_\Omega\\
							&  =(\varphi,\bfu\cdot\bfv+\alpha(\gamma-\gamma_{0}))_\Omega,
\end{align*}
where $A(\gamma)=(\bfu,p)$ and $B(\gamma)=(\bfv,q)$.

\begin{definition}
For the constants $a,b\in\RR$, with $0\leq a<b$, the projection operator $\Pi_{[ a,b]}$ is defined in $L^{2}(\Omega)$ as%
\[
(\forall\varphi\in L^{2}(\Omega))\quad \Pi_{[ a,b]}(\varphi
)(x)=\min\{b,\max\{\varphi(x),a\}\}
\]
\end{definition}

Every local solution $\bar{\gamma}\in\mathcal{A}$ of \eqref{min} that satisfies \eqref{first-order} verifies the identity%
\begin{equation}
\bar{\gamma}=\Pi_{[ a,b]}\left(  \gamma_{0}+\dfrac{1}{\alpha
}\boldsymbol{\bar{u}}\cdot\boldsymbol{\bar{v}}\right), \label{opt-cond}
\end{equation}
where $A(\bar{\gamma})=(\boldsymbol{\bar{u}},\bar{p})$ and $B(\bar{\gamma})=(\boldsymbol{\bar{v}},\bar{q})$.

\begin{proposition}
Let $\varphi_{1},\varphi_{2}\in L^{\infty}(\Omega)$ and $\gamma\in\mathcal{A}$ such that $A^{\prime}(\gamma)[\varphi_{1}]=(\bfu_{1}^{\prime},p_{1}^{\prime})$, $A^{\prime}(\gamma)[\varphi_{2}]=(\bfu_{2}^{\prime},p_{2}^{\prime})$, $A^{\prime\prime}(\gamma)[\varphi_{1},\varphi_{2}]=(\bfu^{\prime\prime},p^{\prime\prime})$. Then, there exists a constant $C>0$ such that

\begin{enumerate}
\item $\Vert\bfu\Vert_{2,\Omega}\leq C\Vert\bff\Vert_{0,\Omega}$.

\item $\Vert\bfu_{j}^{\prime}\Vert_{0,\Omega}\leq C\vert\bfu\vert_{1,\Omega}\Vert\varphi_{j}\Vert_{0,\Omega}$ for $j\in\{1,2\}$.

\item $\Vert\bfu_{j}^{\prime}\Vert_{2,\Omega}\leq C\Vert\varphi_{j}\Vert_{0,\Omega}$ for $j\in\{1,2\}$.

\item $\Vert\bfu^{\prime\prime}\Vert_{0,\Omega}\leq C\vert\bfu\vert_{1,\Omega}\Vert\varphi_{1}\Vert_{0,\Omega}\Vert\varphi_{2}\Vert_{0,\Omega}$.

\item $\Vert\bfu^{\prime\prime}\Vert_{2,\Omega}\leq C\Vert\varphi_{1}\Vert_{0,\Omega}\Vert\varphi_{2}\Vert_{0,\Omega}$.
\end{enumerate}
\end{proposition}

The following result gives us important characteristics about the smoothness of the local solutions of \eqref{min}.
\begin{theorem}\label{H2}
Let $\bar{\gamma}\in\mathcal{A}$ a local solution of \eqref{min} for $\gamma_{0}\in H^{1}\left(  \Omega\right)  $. Then, $\gamma\in H^{1}(\Omega)$ and there exists a positive constant $C>0$ such that
\[
\Vert\bar{\gamma}\Vert_{1,\Omega}\leq\Vert\gamma_{0}\Vert_{1,\Omega}+C\Vert\bff\Vert_{0,\Omega}(\Vert\bff\Vert_{0,\Omega}+\Vert\bfu_{0}\Vert_{0,\Omega}).
\]

\end{theorem}

\begin{proof}
From Lemma \ref{h2}, if $A(\bar{\gamma})=(\boldsymbol{\bar{u}},\bar{p})$ and $B(\bar{\gamma})=(\boldsymbol{\bar{v}},\bar{q})$, then $\boldsymbol{\bar{u}},\boldsymbol{\bar{v}}\in\bfH^{2}(\Omega)$ and $\boldsymbol{\bar {u}},\boldsymbol{\bar{v}}\in\bfL^{\infty}(\Omega)$ by Sobolev Embedding Theorem. Since $\bfH^{2}(\Omega)$ is a Banach algebra for $d\in\{2,3\}$, then $\boldsymbol{\bar{u}}\cdot\boldsymbol{\bar{v}}\in H^{2}(\Omega)$ allow us to deduce that $\nabla(\boldsymbol{\bar{u}}\cdot\boldsymbol{\bar{v}})\in\bfL^{2}(\Omega)$. Since $\Pi_{[a,b]}:H^{1}(\Omega)\rightarrow H^{1}(\Omega)$ is continuous, we can conclude that $\bar{\gamma}\in H^{1}(\Omega)$. The bound is a direct consequence of Lemma \ref{H2}.
\end{proof}

The following results allow us to characterize the second-order optimality condition. The first step is to obtain an expression for the second derivative of $J$.

\begin{proposition}
Let $\varphi\in L^{\infty}(\Omega)$ and $\gamma\in\mathcal{A}$. If $B(\gamma)=(\bfv,q)$ and $A^{\prime}(\gamma)[\varphi]=(\bfu^{\prime},p^{\prime})$, then
\[
J^{\prime\prime}(\gamma)[\varphi,\varphi]=\Vert\bfu^{\prime}\Vert_{0,\omega}^{2}-2(\varphi\bfu^{\prime}+(\nabla\bfu^{\prime})\bfu^{\prime},\bfv)_\Omega+\alpha\Vert\varphi\Vert_{0,\Omega}^{2}.
\]
Furthermore, there exists a positive constant $C>0$ such that for all $\gamma_{1},\gamma_{2}\in\mathcal{A}$%
\[
\vert J^{\prime\prime}(\gamma_{1})[\varphi,\varphi]-J^{\prime\prime}(\gamma_{2})[\varphi,\varphi]\vert\leq C\left\Vert \gamma_{1}-\gamma_{2}\right\Vert_{0,\Omega}\left\Vert \varphi\right\Vert _{0,\Omega}^{2}.
\]
\end{proposition}

\begin{proof}
See Section 6 in \cite{ABO23}.
\end{proof}

\begin{remark}\label{assumption-1}
In that follows, we assume that there exist a positive constant $\delta_{1}>0$ such that for all $\varphi\in L^{\infty}(\Omega)$ and for every local solution $\bar{\gamma}$ that fulfills \eqref{first-order}
\begin{align}
J^{\prime\prime}(\bar{\gamma})[\varphi,\varphi]\geq\delta_{1}\Vert\varphi\Vert_{0,\Omega}^{2}. \label{second-order}
\end{align}
Indeed, this assumption can be fulfilled when $\Vert\boldsymbol{\bar{u}}-\bfu_{0}\Vert_{0,\omega}$ is small enough or $\alpha$ is large enough, since
\begin{align*}
J^{\prime\prime}(\bar{\gamma})[\varphi,\varphi]  &  =\Vert\boldsymbol{\bar{u}}^{\prime}\Vert_{0,\omega}^{2}-(\boldsymbol{\bar{u}}-\bfu_{0},\boldsymbol{\bar{u}}^{\prime\prime})_{\omega}+\alpha\Vert\varphi\Vert_{0,\Omega}^{2}\\
																								 &  \geq\alpha\Vert\varphi\Vert_{0,\Omega}^{2}-\Vert\boldsymbol{\bar{u}}-\bfu_{0}\Vert_{0,\omega}\Vert\boldsymbol{\bar{u}}^{\prime\prime}\Vert_{0,\Omega}\\
																								 &  \geq\alpha\Vert\varphi\Vert_{0,\Omega}^{2}-C\Vert\boldsymbol{\bar{u}}-\bfu_{0}\Vert_{0,\omega}(1+\Vert\bar{\gamma}\Vert_{0,\Omega})\Vert\varphi\Vert_{0,\Omega}^{2}\\
																								 &  \geq(\alpha-C\Vert\boldsymbol{\bar{u}}-\bfu_{0}\Vert_{0,\omega}(1+\Vert\bar{\gamma}\Vert_{0,\Omega}))\Vert\varphi\Vert_{0,\Omega}^{2}.
\end{align*}
\end{remark}

\begin{lemma}\label{ellipticity}
Let $\bar{\gamma}\in\mathcal{A}$ a local solution of \eqref{min} that verifies the Assumption \eqref{second-order}. There exists a constant $\varepsilon>0$ such that for all $\varphi\in L^{\infty}(\Omega)$ and all $\gamma\in\mathcal{A}_{\varepsilon}(\bar{\gamma})$
\[
J^{\prime\prime}(\gamma)[\varphi,\varphi]\geq\dfrac{\delta}{2}\Vert
\varphi\Vert_{0,\Omega}^{2}.
\]
\end{lemma}
\begin{proof}
It is a straightforward consequence of Remark \ref{assumption-1}.
\end{proof}
Finally, we present the second-order optimality condition.
\begin{theorem}
Let $\bar{\gamma}\in\mathcal{A}$ a local solution of \eqref{min} that verifies \eqref{first-order} and the Assumption \eqref{second-order}. There exist positive constants $\sigma$ and $\varepsilon$ such that for all $\gamma\in\mathcal{A}_{\varepsilon}(\bar{\gamma})$.
\[
J(\gamma)\geq J(\bar{\gamma})+\sigma\Vert\gamma-\bar{\gamma}\Vert_{0,\Omega
}^{2}.
\]
\end{theorem}
\begin{proof}
See Theorem 25 in \cite{ABO23}.
\end{proof}

\section{Discretization by finite elements of the states and adjoint equations}\label{S3}
Consider $\Omega\subseteq\RR^{n}$ as a bounded polygonal domain. Let $\left\{  \mathcal{T}_{h}\right\}  _{h>0}$ be a shape-regular family of triangulations of $\overline{\Omega}$ composed of triangles (if $n=2$) or tetrahedron (if $n=3$), with $h:=\max\left\{  h_{T}\mid T\in\mathcal{T}_{h}\right\}  $, where $h_{T}:=\operatorname{diam}(T)$ is the diameter of $T\in\mathcal{T}_{h}$. We define $\mathcal{E}_{h}$ as the set of all edges (faces) of $\mathcal{T}_{h}$.

\begin{definition}\label{discrete-spaces}
We denote $H^{h}$ the continuous Lagrange finite element vector space with degree $2$ on $\overline{\Omega}$, i.e.,%
\[
H^{h}:=\{\bfv^{h}\in\boldsymbol{C}(\overline{\Omega})  \mid(\forall T\in\mathcal{T}_{h})  \quad\bfv^{h}\vert_{T}\in\PP_{2}(T)^{n}\} \cap H,
\]
where $\PP_{2}(T)$ is the space of polynomials of total degree at most $2$ defined on $T$. Analogously, we denote $Q^{h}$ the continuous Lagrange finite element space with degree $1$ on $\overline{\Omega}$, that is,
\[
Q^{h}:=\{q^{h}\in C(\overline{\Omega})  \mid(\forall T\in\mathcal{T}_{h})  \quad q^{h}\vert_{T}\in\PP_{1}(T)\} \cap Q,
\]
where $\PP_{1}(T)$ is the space of polynomials of total degree at most $1$ defined on $T$. We denote the Taylor-Hood finite element as $V^{h}=H^{h}\times Q^{h}$.
\end{definition}

\begin{theorem}\label{taylor-hood}
Let $\gamma\in\mathcal{A}$, $\bfu_{0}\in\bfL^{2}(\omega)$, $\gamma_{0}\in L^{2}(\Omega)$ and $\bff\in\bfL^{2} (\Omega)$ such that verifies the hypothesis of Theorem \ref{f-bound}. If the solutions $(\bfu,p)$ for \eqref{VF} and $(\bfv,q)$ for \eqref{AVF} are smooth enough, i.e. $\bfu,\bfv\in\bfH^{3}(\Omega)\cap H$ and $p,q\in H^{2}(\Omega)\cap Q$, then the solution of the discrete variational formulations given by

\textit{Find $(\bfu^h,p^h),(\bfv^h,q^h)\in V^{h}$ such that for all $(\bfw^{h},r^{h})\in V^{h}$}
\begin{align}
a(\bfu^{h},\bfw^{h})-b(\bfw^{h},p^{h})+b(\bfu^{h},r^{h})+c(\bfu^{h},\bfu^{h},\bfw^{h})+(\gamma\bfu^{h},\bfv^{h})_\Omega  & =(\bff,\bfw^{h})_\Omega\label{VFh}\\
a(\bfv^{h},\bfw^{h})-b(\boldsymbol{r}^{h},q^{h})+b(\bfv^{h},r^{h})+c(\bfu^{h},\bfv^{h},\bfw^{h})+c(\bfv^{h},\bfu^{h},\bfw^{h})+(\gamma\bfv^{h},\bfw^{h})_\Omega  &  =(\bfu^{h}-\bfu_{0},\bfw^{h})_{\omega} \label{AVFh}
\end{align}
is unique and the following error estimates hold
\begin{align*}
\Vert (\bfu-\bfu^{h},p-p^{h})\Vert    &  \leq Ch^2 (\Vert\bfu\Vert_{3,\Omega}+\Vert p\Vert_{2,\Omega})\\
\Vert \bfu-\bfu^{h}\Vert _{0,\Omega}  &  \leq Ch^{3}(\Vert\bfu\Vert_{3,\Omega}+\Vert p\Vert_{2,\Omega})\\
\Vert (\bfv-\bfv^{h},q-q^{h})\Vert    &  \leq Ch^2 (\Vert\bfv\Vert_{3,\Omega}+\Vert q\Vert_{2,\Omega})\\
\Vert \bfv-\bfv^{h}\Vert _{0,\Omega}  &  \leq Ch^{3} (\Vert\bfv\Vert_{3,\Omega}+\Vert q\Vert_{2,\Omega}).
\end{align*}
\end{theorem}

\begin{proof}
See Sections 5.2 and 6.2 in \cite{J16}.
\end{proof}

\begin{definition}
We denote as $A_{h}:\mathcal{A}\rightarrow V^h$ and $B_{h}:\mathcal{A}\rightarrow V^h$ are the discretized control-to-state and control-to-adjoint-state operators such that $A_{h}(\gamma)=(\bfu^{h},p^{h})$ and $B_{h}(\gamma)=(\bfv^{h},q^{h})$ are the solutions of \eqref{VFh} and \eqref{AVFh}, respectively. Analogously, we define the functional $J_{h}:\mathcal{A}\rightarrow\RR$ given by $J_{h}(\gamma)=\dfrac{1}{2}\Vert\bfu_{h}-\bfu_{0}\Vert_{0,\omega}^{2}+\dfrac{\alpha}{2}\Vert\gamma-\gamma_{0}\Vert_{0,\Omega}^{2}$, where $A_{h}(\gamma)=(\bfu^{h},p^{h})$.
\end{definition}

\begin{remark}
The Fr\'{e}chet-differentiability of $A_{h}$, $B_{h}$ and $J_{h}$ is inherited from $A$, $B$ and $J$, respectively.
\end{remark}

\begin{lemma}
Let $\gamma\in\mathcal{A}$, $A_{h}(\gamma)=(\bfu^{h},p^{h})$ and $B_{h}(\gamma)=(\bfv^{h},q^{h})$. Then, $\Vert\bfu^{h}\Vert_{0,\infty,\Omega}\leq C\Vert\bff\Vert_{0,\Omega}$ and $\Vert\bfu^{h}\Vert_{0,\infty,\Omega}\leq C(\Vert\bfu_{0}\Vert_{0,\omega}+\Vert\bff\Vert_{0,\Omega})$.
\end{lemma}

\begin{proof}
By Sobolev Embedding Theorem, Theorem 3.24 from \cite{EG04} and Theorem \ref{f-bound}, we have%
\[
\Vert\bfu^{h}\Vert_{0,\infty,\Omega}\leq C\Vert\bfu%
^{h}\Vert_{1,4,\Omega}\leq C\vert\bfu\vert_{1,4,\Omega}\leq
C\vert\bfu\vert_{2,\Omega}\leq C\Vert\bff\Vert_{0,\Omega}%
\]
The second estimate is obtained by the same way.
\end{proof}

\begin{lemma}
Let $\gamma\in\mathcal{A}$ and $\varphi\in L^{\infty}(\Omega)$. Then,

\begin{enumerate}
\item $\vert A^{\prime}(\gamma)[\varphi]-A_{h}^{\prime}(\gamma)[\varphi]\vert\leq Ch\Vert\varphi\Vert_{0,\Omega}(\Vert\bff\Vert_{0,\Omega}+\Vert\bff\Vert_{0,\Omega}^{2})$.

\item $\vert A^{\prime\prime}(\gamma)[\varphi,\varphi]-A_{h}^{\prime\prime}(\gamma)[\varphi,\varphi]\vert\leq Ch^{2}\Vert\varphi\Vert_{0,\Omega}^{2}(\Vert\bff\Vert_{0,\Omega}+\Vert\bff\Vert_{0,\Omega}^{2})^{2}$.
\end{enumerate}
\end{lemma}

\begin{proof}
If $A(\gamma)=(\bfu,p)$ and $A_{h}(\gamma)=(\bfu^{h},p^{h})$, we have that $A^{\prime}(\gamma)[\varphi]=(\bfu^{\prime},p^{\prime})$ and $A_{h}^{\prime}(\gamma)[\varphi]=(\bfu_{h}^{\prime},p_{h}^{\prime})$ verifies the identities
\begin{align*}
a(\bfu^{\prime},\bfv)-b(\bfv,p^{\prime})+b(\bfu^{\prime},q)+c(\bfu^{\prime},\bfu,\bfv)+c(\bfu,\bfu^{\prime},\bfv)+(\gamma \bfu^{\prime},\bfv)_\Omega=-(\varphi\bfu,\bfv)_\Omega,\\
a(\bfu_{h}^{\prime},\bfv^{h})-b(\bfv^{h},p_{h}^{\prime})+b(\bfu_{h}^{\prime},q^{h})+c(\bfu_{h}^{\prime},\bfu^{h},\bfv^{h})+c(\bfu^{h},\bfu_{h}^{\prime},\bfv^{h})+(\gamma\bfu_{h}^{\prime},\bfv^{h})_\Omega=-(\varphi\bfu^{h},\bfv^{h})_\Omega,
\end{align*}
for all $(\bfv,q) \in V$ and $(\bfv^h,q^h) \in V$. If $(\boldsymbol{\hat{u}}^{h},\hat{p}^{h})$ are the solutions of%
\[
(\forall(\bfv^{h},q^{h})\in H^{h}\times Q^{h})\quad a(\boldsymbol{\hat{u}}^{h},\bfv^{h})-b(\bfv^{h},\hat{p}^{h})+b(\boldsymbol{\hat{u}}^{h},q^{h})+c(\boldsymbol{\hat{u}}^{h},\bfu,\bfv^{h})+c(\bfu,\boldsymbol{\hat{u}}^{h},\bfv^{h})+(\gamma\boldsymbol{\hat{u}}^{h},\bfv%
^{h})_\Omega=-(\varphi\bfu,\bfv^{h})_\Omega.
\]
Then,%
\begin{align*}
\Vert (\bfu^{\prime}-\boldsymbol{\hat{u}}^{h},p^{\prime}-\hat{p}^{h})\Vert  &  \leq Ch\Vert\varphi\Vert_{0,\Omega}\vert\bfu\vert_{1,\Omega}\\
&  \leq Ch\Vert\varphi\Vert_{0,\Omega}\Vert\bff\Vert_{0,\Omega},\\
\Vert (\bfu_{h}^{\prime}-\boldsymbol{\hat{u}}^{h},p_{h}^{\prime}-\hat{p}^{h})\Vert  &  \leq C\Vert\varphi\Vert_{0,\Omega}\Vert\bfu-\bfu^{h}\Vert_{1,\Omega}\\
&  \leq Ch\Vert\varphi\Vert_{0,\Omega}(\vert\bfu\vert_{2,\Omega}+\Vert p\Vert_{1,\Omega})\\
&  \leq Ch\Vert\varphi\Vert_{0,\Omega}(\Vert\bff\Vert_{0,\Omega}+\Vert\bff\Vert_{0,\Omega}^{2}).
\end{align*}
By Triangle inequality, we conclude
\[
\Vert A^{\prime}(\gamma)[\varphi]-A_{h}^{\prime}(\gamma)[\varphi]\Vert\leq Ch\Vert\varphi\Vert_{0,\Omega}(\Vert\bff\Vert_{0,\Omega}+\Vert\bff\Vert_{0,\Omega}^{2}).
\]
The second estimate is obtained following similar steps.
\end{proof}

\begin{lemma}\label{aux-1}
Let $\gamma,\gamma_{1},\gamma_{2}\in\mathcal{A}$, $\varphi\in L^{\infty}(\Omega)$. Then,
\begin{enumerate}
\item $\left\vert J^{\prime}(\gamma)[\varphi]-J_{h}^{\prime}(\gamma)[\varphi]\right\vert \leq Ch^{2}\Vert\varphi\Vert_{0,\Omega}\Vert\bff\Vert_{0,\Omega}(\Vert\bff\Vert_{0,\Omega}+\Vert\bfu_{0}\Vert_{0,\omega})$.
\item $\left\vert J_{h}^{\prime}(\gamma_{1})[\varphi]-J_{h}^{\prime}(\gamma_2)[\varphi]\right\vert \leq C\Vert\gamma_{1}-\gamma_{2}\Vert_{0,\Omega}\Vert\varphi\Vert_{0,\Omega}$.
\end{enumerate}
\end{lemma}

\begin{proof}
If $A(\gamma)=(\bfu,p)$, $A_{h}(\gamma)=(\bfu^{h},p^{h})$, $B(\gamma)=(\bfv,q)$ and $B_{h}(\gamma)=(\bfv^{h},q^{h})$, then
\begin{align*}
\left\vert J^{\prime}(\gamma)[\varphi]-J_{h}^{\prime}(\gamma)[\varphi]\right\vert  &  =\vert (\varphi,\bfu\cdot\bfv -\bfu^{h}\cdot\bfv^{h})_\Omega \vert \\
&  =\vert (\varphi,(\bfu-\bfu^{h})\cdot\bfv)_\Omega +(\varphi,\bfu^{h}\cdot(\bfv-\bfv^{h}))_\Omega\vert \\
&  \leq\vert (\varphi,(\bfu-\bfu^{h})\cdot \bfv)_\Omega\vert+\vert(\varphi,\bfu^{h}\cdot(\bfv-\bfv^{h}))\vert \\
&  \leq C\Vert\varphi\Vert_{0,\Omega}(\Vert\bfu-\bfu^{h}\Vert_{0,\Omega}\vert\bfv\vert_{0,\infty,\Omega}+\Vert\bfv-\bfv^{h}\Vert_{0,\Omega}\Vert\bfu^{h}\Vert_{0,\infty,\Omega})\\
&  \leq Ch^{2}\Vert\varphi\Vert_{0,\Omega}\Vert\bff\Vert_{0,\Omega}(\Vert\bff\Vert_{0,\Omega}+\vert\bfu\vert_{0,\omega}),
\end{align*}
proving the first estimate. The second estimate is similar to the proof of Lemma \ref{lipschitz}.
\end{proof}

\begin{lemma}\label{ellipticity-h}
Let $\bar{\gamma}\in\mathcal{A}$ a local solution of \eqref{min} that verifies the Assumption \eqref{second-order}. There exists $h_0>0$ and a constant $\varepsilon>0$ such that for all $h\in(0,h_0)$, for all $\varphi\in L^{\infty}(\Omega)$ and for all $\gamma\in\mathcal{A}_{\varepsilon}(\bar{\gamma})$
\[
J_{h}^{\prime\prime}(\gamma)[\varphi,\varphi]\geq\dfrac{\delta}{4}\Vert
\varphi\Vert_{0,\Omega}^{2}.
\]
\end{lemma}

\begin{proof}
First, we have $J^{\prime\prime}(\gamma)[\varphi,\varphi]\geq \dfrac{\delta}{2}\Vert\varphi\Vert_{0,\Omega}^{2}$ from Lemma \ref{ellipticity}. Then, from Lemma \ref{aux-1},
\begin{align*}
\vert J_{h}^{\prime\prime}(\gamma)[\varphi,\varphi]-J^{\prime\prime}(\gamma)[\varphi,\varphi]\vert  &  =\vert\Vert\bfu_{h}^{\prime}\Vert_{0,\omega}^{2}-(\bfu^{h}-\bfu_{0},\bfu_{h}^{\prime\prime})_{\omega}-(\Vert\bfu^{\prime}\Vert_{0,\omega}^{2}-(\bfu-\bfu_{0},\boldsymbol{\bar{u}}^{\prime\prime})_{\omega})\vert\\
&  \leq\vert\Vert\bfu_{h}^{\prime}\Vert_{0,\omega}^{2}-\Vert\bfu^{\prime}\Vert_{0,\omega}^{2}\vert+\vert(\bfu-\bfu_{0},\bfu^{\prime\prime})_{\omega}-(\bfu^{h}-\bfu_{0},\bfu_{h}^{\prime\prime})_{\omega}\vert\\
&  \leq\vert\Vert\bfu_{h}^{\prime}\Vert_{0,\omega}^{2}-\Vert\bfu^{\prime}\Vert_{0,\omega}^{2}\vert+\vert(\bfu-\bfu^{h},\bfu^{\prime\prime})_{\omega}\vert+\vert(\bfu^{h}-\bfu_{0},\bfu_{h}^{\prime\prime})_{\omega}\vert\\
&  \leq Ch\Vert\varphi\Vert_{0,\Omega}^{2}(\Vert\bff\Vert_{0,\Omega}+\Vert\bff\Vert_{0,\Omega}^{2})^{2}.
\end{align*}
Applying Triangle inequality,
\begin{align*}
J_{h}^{\prime\prime}(\gamma)[\varphi,\varphi]  &  =J^{\prime\prime}(\gamma)[\varphi,\varphi]+J_{h}^{\prime\prime}(\gamma)[\varphi,\varphi]-J^{\prime\prime}(\gamma)[\varphi,\varphi]\\
&  \geq\left(  \dfrac{\delta}{2}-Ch(\Vert\bff\Vert_{0,\Omega}+\Vert\bff\Vert_{0,\Omega}^{2})^{2}\right)  \Vert\varphi\Vert_{0,\Omega}^{2}.
\end{align*}
Choosing $h_0=\dfrac{\delta}{4(\Vert\bff\Vert_{0,\Omega} +\Vert\bff\Vert_{0,\Omega}^{2})^{2}}$, we have that $\left(\dfrac{\delta}{2}-Ch(\Vert\bff\Vert_{0,\Omega}+\Vert\bff\Vert_{0,\Omega}^{2})^{2}\right)  \geq\dfrac{\delta}{4}$ for any $h\in(0,h_0)$, concluding that $J_{h}^{\prime\prime}(\gamma)[\varphi,\varphi]\geq\dfrac{\delta}{4}\Vert\varphi\Vert_{0,\Omega}^{2}$.
\end{proof}

\begin{theorem}
The semi-discrete optimal control problem
\begin{align}
\text{minimize\quad}  &  J_{h}(\gamma)=\dfrac{1}{2}\Vert\bfu^{h}-\bfu_{0}\Vert_{0,\omega}^{2}+\dfrac{\alpha}{2}\Vert\gamma-\gamma_{0}\Vert_{0,\Omega}^{2} \label{semi-h} \\
\text{subject to\quad}  &  (\forall(\bfw^{h},r^{h})\in V^h)\quad a(\bfu^{h},\bfw^{h})-b(\bfw^{h},p^{h})+b(\bfu^{h},r^{h})+c(\bfu^{h},\bfu^{h},\bfw^{h})+(\gamma\bfu^{h},\bfv^{h})=(\bff,\bfw^{h})\nonumber\\
&  \gamma\in\mathcal{A}\nonumber
\end{align}
admits at least one global solution $(\boldsymbol{\bar{u}}^{h},\bar{p}^{h},\bar{\gamma})\in H^{h}\times Q^{h}\times\mathcal{A}$.
\end{theorem}
\begin{proof}
The proof is similar to the one for Theorem \ref{existence}. We omit the details.
\end{proof}

\subsection{Discretization of the optimal control problem}

\begin{definition}
We denote
\begin{align*}
G_{0}^{h}  &  :=\left\{  \gamma^{h}\in L^{2}(\Omega)\mid\left(  \forall T\in\mathcal{T}_{h}\right)  (\exists a_{T}\in\RR)\quad\gamma^{h}\vert_{T}=a_{T}\right\}, \\
G_{1}^{h}  &  :=\left\{  \gamma^{h}\in C\left(  \overline{\Omega}\right) \mid\left(  \forall T\in\mathcal{T}_{h}\right)  \quad\gamma^{h}\vert_{T}\in\PP_{1}(T)\right\}.
\end{align*}

\end{definition}

The spaces $G_{0}^{h}$ and $G_{1}^{h}$ are some possible finite element spaces for the optimal control, corresponding to discontinuous and continuous Galerkin finite element spaces with the lowest polynomial degree. In that case, $\mathcal{A}$ can be changed by $\mathcal{A}^{h}=\mathcal{A}\cap G_{0}^{h}$ or $\mathcal{A}^{h}=\mathcal{A}\cap G_{1}^{h}$ in the semi-discrete optimal control problem obtaining a fully discrete problem. This problem also has at least one global solution $(\boldsymbol{\bar{u}}^{h},\bar{p}^{h},\bar{\gamma}^{h})\in H^{h}\times Q^{h}\times\mathcal{A}^{h}$. 

In the first two subsections, we consider the cases $\mathcal{A}^{h}=\mathcal{A}\cap G_{0}^{h}$ or $\mathcal{A}^{h}=\mathcal{A}\cap G_{1}^{h}$. In both cases, the optimal control $\bar{\gamma}^h$ must fulfull a discrete version of \eqref{first-order} given by
\begin{equation}
(\forall \varphi^h \in \mathcal{A}^h)\quad (\varphi^h-\bar{\gamma}^h,\boldsymbol{u}^h\cdot\boldsymbol{v}^h+\alpha(\bar{\gamma}^h-\gamma_{0}))_\Omega \geq 0. \label{first-order-h}
\end{equation}
In the final subsection, we also analyze one alternative that consists of a semi-discrete optimal control problem, where only the states and adjoints are discretized. The optimal control must also fulfills the optimality conditions obtained the previous sections. Then, every local solution $\bar{\gamma}^{h}\in\mathcal{A}$ must verify the identity
\begin{equation}
\bar{\gamma}^{h}=\Pi_{[ a,b]}\left(  \gamma_{0}+\dfrac{1}{\alpha
}\boldsymbol{\bar{u}}^{h}\cdot\boldsymbol{\bar{v}}^{h}\right). \label{gamma-h}
\end{equation}

\begin{definition}
We define the discrete optimal control problem as
\begin{align}
\text{minimize\quad}  &  J_{h}(\gamma^h)=\dfrac{1}{2}\Vert\bfu^{h}-\bfu_{0}\Vert_{0,\omega}^{2}+\dfrac{\alpha}{2}\Vert\gamma^{h}-\gamma_{0}\Vert_{0,\Omega}^{2} \label{min-h}\\
\text{subject to\quad}  &  (\forall(\bfw^{h},r^{h})\in V^{h})\quad a(\bfu^{h},\bfw^{h})-b(\bfw^{h},p^{h})+b(\bfu^{h},r^{h})+c(\bfu^{h},\bfu^{h},\bfw^{h})+(\gamma^{h}\bfu^{h},\bfv^{h})_\Omega=(\bff,\bfw^{h})_\Omega \nonumber\\
&  \gamma^{h}\in\mathcal{A}^{h}.\nonumber
\end{align}
\end{definition}

In order to obtain our a priori error estimates, we need to define a local auxiliary optimal control problem.
\begin{definition}
Consider the constants $\varepsilon>0$ and $h>0$. Given $\gamma^\ast \in\mathcal{A}$, we denote%
\[
\mathcal{A}_{\varepsilon}^{h}(\gamma^\ast)=\{\gamma^{h}\in\mathcal{A}%
^{h}\text{ }\mid\text{ }\Vert\gamma^{h}-\gamma^\ast\Vert _{0,\Omega} \leq\varepsilon\}.
\]
Then, for $\bar{\gamma}\in\mathcal{A}$ a local solution of \eqref{min}, we define the auxiliary optimal control problem as
\begin{align}
\text{minimize\quad}  &  J_{h}(\gamma^{h})=\dfrac{1}{2}\Vert\bfu^{h}-\bfu_{0}\Vert_{0,\omega}^{2}+\dfrac{\alpha}{2}\Vert\gamma^{h}-\gamma_{0}\Vert_{0,\Omega}^{2} \label{aux-oc}\\
\text{subject to\quad}  &  (\forall(\bfw^{h},r^{h})\in V^{h})\quad a(\bfu^{h},\bfw^{h})-b(\bfw^{h},p^{h})+b(\bfu^{h},r^{h})+c(\bfu%
^{h},\bfu^{h},\bfw^{h})+(\gamma^{h}\bfu^{h},\bfv^{h})_\Omega =(\bff,\bfw^{h})_\Omega \nonumber\\
&  \gamma^{h}\in\mathcal{A}_{\varepsilon}^{h}(\bar{\gamma}).\nonumber
\end{align}
\end{definition}

Under some hypotheses, this auxiliary control problem has an unique solution.

\begin{definition}
We denote by $\mathcal{P}_{0}^{h}:L^{2}(\Omega)\rightarrow G_{0}^{h}$ and $\mathcal{I}_{1}^{h}:C(\Omega)\rightarrow G_{1}^{h}$ the orthogonal projection to $G_{0}^{h}$ with respect to the $L^{2}(\Omega)$ inner product and the Lagrange interpolator for $G_{1}^{h}$, respectively.
\end{definition}

\begin{proposition}\label{interpolator}
Let $\varphi\in H^{1}(\Omega)$. Then,
\begin{enumerate}
\item $\Vert\varphi-\mathcal{P}_{0}^{h}(\varphi)\Vert_{0,\Omega}\leq Ch\Vert\varphi\Vert_{1,\Omega}$. Furthermore, if $\varphi\in\mathcal{A}$, then $\mathcal{P}_{0}^{h}(\varphi)\in\mathcal{A}\cap G_{0}^{h}$.

\item If $\varphi\in H^{2}(\Omega)$, then $\Vert\varphi-\mathcal{I}_{1}^{h}(\varphi)\Vert_{1,\Omega}\leq Ch^{2}\Vert\varphi\Vert_{2,\Omega}$. Furthermore, if $\varphi\in\mathcal{A}$, then $\mathcal{I}_{1}^{h}(\varphi)\in\mathcal{A}\cap G_{1}^{h}$.
\end{enumerate}
\end{proposition}

\begin{proof}
See Proposition 1.134 and Corollary 1.109 in \cite{EG04}.
\end{proof}

\begin{lemma}\label{existence-h}
Consider $\mathcal{A}^{h}=\mathcal{A}\cap G_{0}^{h}$ or $\mathcal{A}^{h}=\mathcal{A}\cap G_{1}^{h}$. For all $\varepsilon>0$, there exists $h_0>0$ such that the auxiliary optimal control problem has a solution for all $h\in(0,h_0)$.
\end{lemma}

\begin{proof}
Taking $\bar{\gamma}^{h}=\mathcal{P}_{0}^{h}(\bar{\gamma})$ if $\mathcal{A}^{h}=\mathcal{A}\cap G_{0}^{h}$, or $\bar{\gamma}^{h} =\mathcal{I}_{1}^{h}(\bar{\gamma})$ if $\mathcal{A}^{h}=\mathcal{A}\cap G_{1}^{h}$, there exists $h_0 >0$ such that $\Vert\bar{\gamma}-\bar{\gamma}^{h}\Vert_{0,\Omega}\leq\varepsilon$. Then, $\mathcal{A}_{\varepsilon}^{h}(\bar{\gamma})$ is non empty. We can reply the same techniques as in \cite{ABO23} to prove the existence of a control.
\end{proof}

\begin{remark}\label{local-h}
One of the straightforward consequences of the proof of Lemma \ref{existence-h} is that it describes a way to obtain a local solution for the discrete problem. Defining 
\begin{equation}
\gamma^{\ast} = \Pi_{[ a,b]}\left(  \gamma_{0}+\dfrac{1}{\alpha
}\boldsymbol{\bar{u}}^{h}\cdot\boldsymbol{\bar{v}}^{h}\right),
\end{equation}
we have that $\gamma^{h}=\mathcal{P}_{0}^{h}(\gamma^\ast)$ and $\gamma^{h}=\mathcal{I}_{1}^{h}(\gamma^\ast)$ are local solutions for $\mathcal{A}^{h}=\mathcal{A}\cap G_{0}^{h}$ and $\mathcal{A}^{h}=\mathcal{A}\cap G_{1}^{h}$, respectively. Furthermore, in the case $\mathcal{A}^{h}=\mathcal{A}\cap G_{0}^{h}$, \eqref{first-order-h} can be reduced to the following variational formulation
\begin{equation}
(\forall \varphi^h \in G_0 ^{h}) \quad (\gamma^h - \gamma^\ast,\varphi^h)_\Omega=0. \label{first-order-p0}
\end{equation}
However, we cannot reply this same reasoning when $\mathcal{A}^{h}=\mathcal{A}\cap G_{1}^{h}$ since $\mathcal{P}_{0}^{h}(\varphi)$ could not belong to $\mathcal{A}$ for all $\varphi\in\mathcal{A}$.
\end{remark}

\begin{lemma}\label{taylor}
Consider $\varepsilon>0$ small enough such that $J_{h}^{\prime\prime}$ is coercive as in Lemma \ref{ellipticity-h} for all $\gamma^{h}\in\mathcal{A}_{\varepsilon}^{h}(\bar{\gamma})$. Then, the auxiliary optimal control problem \eqref{aux-oc} has an unique solution for $h>0$ small enough.
\end{lemma}

\begin{proof}
Let $\gamma_{1}^{h},\gamma_{2}^{h}\in\mathcal{A}_{\varepsilon}^{h} (\bar{\gamma})$ two solutions of \eqref{aux-oc}. Then, by Taylor Theorem, there exists $t\in [0,1]$ such that
\[
J_{h}(\gamma_{1}^{h})=J_{h}(\gamma_{2}^{h})+J_{h}^{\prime}(\gamma_{2}^{h})[\gamma_{1}^{h}-\gamma_{2}^{h}]+\dfrac{1}{2}J^{\prime\prime}(t\gamma_{1}^{h}+(1-t)\gamma_{2}^{h})[\gamma_{1}^{h}-\gamma_{2}^{h},\gamma_{1}^{h}-\gamma_{2}^{h}].
\]
Since $J_{h}(\gamma_{1}^{h})=J_{h}(\gamma_{2}^{h})$ and $J_{h}^{\prime}(\gamma_{2}^{h})[\gamma_{1}^{h}-\gamma_{2}^{h}]\geq0$, we apply Lemma \ref{ellipticity-h}. Thus,%
\[
\dfrac{\delta}{8}\Vert\gamma_{1}^{h}-\gamma_{2}^{h}\Vert_{0,\Omega}^{2}\leq\dfrac{1}{2}J^{\prime\prime}(t\gamma_{1}^{h}+(1-t)\gamma_{2}^{h})[\gamma_{1}^{h}-\gamma_{2}^{h},\gamma_{1}^{h}-\gamma_{2}^{h}]\leq 0,
\]
proving that $\gamma_{1}^{h}=\gamma_{2}^{h}$.
\end{proof}

\subsection{Case with a discontinuous Galerkin discrete control}

In this subsection, we deduce an a priori error estimate for the case $\mathcal{A}^{h}=\mathcal{A}\cap G_{0}^{h}$ also called $\PP_0$ case. First, we present an estimate for the solution of the auxiliary optimal control problem.

\begin{theorem}\label{P0}
Let $\bar{\gamma}\in\mathcal{A}$ a local solution of \eqref{min} that verifies the optimality conditions \eqref{first-order} and \eqref{second-order}. There exist $\varepsilon>0$ and $h>0$ such that the auxiliary optimal control problem \eqref{aux-oc} has an unique solution $\gamma_{\varepsilon}^{h}\in\mathcal{A}_{\varepsilon}^{h}(\bar{\gamma})$. Furthermore, there exists $h_0>0$ and a positive constant $C$ such that for all $h\in (0,h_0)$
\[
(\forall h\in(0,h_0))\quad\Vert\bar{\gamma}-\gamma_{\varepsilon}^{h}\Vert_{0,\Omega}\leq\dfrac{C}{\delta^{1/2}}h((\Vert\bfu_{0}\Vert_{0,\omega}+\Vert\bff\Vert_{0,\Omega})^{2}+\alpha\Vert\gamma_{0}\Vert_{1,\Omega}).
\]
\end{theorem}

\begin{proof}
From Lemma \ref{ellipticity}, there exists a constant $\varepsilon>0$ such that for all $\varphi\in L^{\infty}(\Omega)$ and all $\gamma\in\mathcal{A}_{\varepsilon}(\bar{\gamma})$%
\[
J^{\prime\prime}(\gamma)[\varphi,\varphi]\geq\dfrac{\delta}{2}\Vert\varphi\Vert_{0,\Omega}^{2}.
\]
Analogously, there exists $h_0>0$ such that for all $h\in(0,h_0)$, $\varphi\in L^{\infty}(\Omega)$ and $\gamma^{h}\in\mathcal{A}_{\varepsilon
}^{h}(\bar{\gamma})$%
\[
J_{h}^{\prime\prime}(\gamma^{h})[\varphi,\varphi]\geq\dfrac{\delta}{4}%
\Vert\varphi\Vert_{0,\Omega}^{2}.
\]
If $\hat{\gamma}_{\varepsilon}^{h}\in\mathcal{A}_{\varepsilon}^{h}(\bar {\gamma})$ is the unique solution of the following optimal control problem%
\begin{align*}
\text{minimize\quad}  &  J(\gamma)=\dfrac{1}{2}\Vert\bfu-\bfu_{0}\Vert_{0,\omega}^{2}+\dfrac{\alpha}{2}\Vert\gamma
-\gamma_{0}\Vert_{0,\Omega}^{2}\\
\text{subject to\quad}  &  (\forall(\bfw,r)\in V)\quad a(\bfu,\bfw)-b(\bfw,p)+b(\bfu,r)+c(\bfu,\bfu,\bfw)+(\gamma^{h}\bfu,\bfv)_\Omega=(\bff,\bfw)_\Omega\\
&  \gamma^{h}\in\mathcal{A}_{\varepsilon}^{h}(\bar{\gamma}),
\end{align*}
we have
\[
\Vert\bar{\gamma}-\gamma_{\varepsilon}^{h}\Vert_{0,\Omega}\leq\Vert\bar{\gamma}-\hat{\gamma}_{\varepsilon}^{h}\Vert_{0,\Omega}+\Vert\hat{\gamma}_{\varepsilon}^{h}-\gamma_{\varepsilon}^{h}\Vert_{0,\Omega}.
\]
Let $e_{1}=\bar{\gamma}-\hat{\gamma}_{\varepsilon}^{h}$. By Taylor Theorem, there exists $t\in[0,1]$ and $\xi_{t}=t\bar{\gamma}+(1-t)\hat{\gamma }_{\varepsilon}^{h}\in\mathcal{A}_{\varepsilon}(\bar{\gamma})$ such that
\begin{align*}
J^{\prime\prime}(\xi_{t})\left[  e_{1},e_{1}\right]   &  =J^{\prime}(\bar{\gamma})[e_{1}]-J^{\prime}(\hat{\gamma}_{\varepsilon}^{h})[e_{1}]\\ &  =J^{\prime}(\bar{\gamma})[\bar{\gamma}-\hat{\gamma}_{\varepsilon}^{h}]-J^{\prime}(\hat{\gamma}_{\varepsilon}^{h})[\bar{\gamma}-\mathcal{P}_{0}\bar{\gamma}]+J^{\prime}(\hat{\gamma}_{\varepsilon}^{h})[\hat{\gamma}_{\varepsilon}^{h}-\mathcal{P}_{0}\bar{\gamma}].
\end{align*}
Since $\bar{\gamma}$ and $\hat{\gamma}_{\varepsilon}^{h}$ verify the respective optimality conditions, we have $J^{\prime}(\bar{\gamma})[
\bar{\gamma}-\hat{\gamma}_{\varepsilon}^{h}]  \leq0$ and $J^{\prime}(\hat{\gamma}_{\varepsilon}^{h})[\hat{\gamma}_{\varepsilon}^{h}-\mathcal{P}_{0}\bar{\gamma}]\leq0$. Then,
\begin{align*}
\dfrac{\delta}{2}\Vert e_{1}\Vert_{0,\Omega}^{2}  &  \leq J^{\prime\prime}(\xi_{t})\left[  e_{1},e_{1}\right] \\
&  \leq-J^{\prime}(\hat{\gamma}_{\varepsilon}^{h})[\bar{\gamma}-\mathcal{P}_{0}\bar{\gamma}]=-(\alpha(\hat{\gamma}_{\varepsilon}^{h}-\gamma_{0})-\bfu_{\varepsilon}^{h}\cdot\bfv_{\varepsilon}^{h},\bar{\gamma}-\mathcal{P}_{0}\bar{\gamma})_\Omega\\
&  \leq-\alpha(\hat{\gamma}_{\varepsilon}^{h}-\gamma_{0},\bar{\gamma}-\mathcal{P}_{0}\bar{\gamma})_\Omega +(\bfu_{\varepsilon}^{h}\cdot\bfv_{\varepsilon}^{h},\bar{\gamma}-\mathcal{P}_{0}\bar{\gamma})_\Omega,
\end{align*}
where $A(\hat{\gamma}_{\varepsilon}^{h})=(\bfu_{\varepsilon},p_{\varepsilon})$ and $B(\hat{\gamma}_{\varepsilon})=(\bfv_{\varepsilon},q_{\varepsilon})$. Since
\[
(\mathcal{P}_{0}(\bfu_{\varepsilon}\cdot\bfv_{\varepsilon}),\bar{\gamma}-\mathcal{P}_{0}\bar{\gamma}) _\Omega =(\mathcal{P}_{0}(\gamma_{0}),\bar{\gamma}-\mathcal{P}_{0}\bar{\gamma}) _\Omega =(\hat{\gamma}_{\varepsilon}^{h},\bar{\gamma}-\mathcal{P}_{0}\bar{\gamma})_\Omega =0,
\]
we have
\begin{align*}
-\alpha(\hat{\gamma}_{\varepsilon}^{h}-\gamma_{0},\bar{\gamma}-\mathcal{P}_{0}(\bar{\gamma}))_\Omega  &  =\alpha(\mathcal{P}_{0}(\gamma_{0})-\gamma_{0},\bar{\gamma}-\mathcal{P}_{0}\bar{\gamma}) _\Omega\\
&  \leq\dfrac{\alpha^{2}}{2}\Vert\gamma_{0}-\mathcal{P}_{0}(\gamma_{0})\Vert_{0,\Omega}^{2}+\dfrac{1}{2}\Vert\bar{\gamma}-\mathcal{P}_{0}\bar{\gamma}\Vert_{0,\Omega}^{2},
\end{align*}
and
\begin{align*}
(\bfu_{\varepsilon}\cdot\bfv_{\varepsilon},\bar{\gamma}-\mathcal{P}_{0}\bar{\gamma})_\Omega  &  =(\bfu_{\varepsilon}\cdot\bfv_{\varepsilon}-\mathcal{P}_{0}(\bfu_{\varepsilon}\cdot\bfv_{\varepsilon}),\bar{\gamma}-\mathcal{P}_{0}\bar{\gamma})_\Omega\\
&  \leq\dfrac{1}{2}\Vert\bfu_{\varepsilon}\cdot\bfv_{\varepsilon}-\mathcal{P}_{0}(\bfu_{\varepsilon}\cdot\bfv_{\varepsilon})\Vert_{0,\Omega}^{2}+\dfrac{1}{2}\Vert\bar{\gamma}-\mathcal{P}_{0}\bar{\gamma}\Vert_{0,\Omega}^{2}.
\end{align*}
Then,%
\begin{align*}
\dfrac{\delta}{2}\Vert e_{1}\Vert_{0,\Omega}^{2}  &  \leq\dfrac{1}{2}\Vert\bfu_{\varepsilon}\cdot\bfv_{\varepsilon}-\mathcal{P}_{0}(\bfu_{\varepsilon}\cdot\bfv_{\varepsilon})\Vert_{0,\Omega}^{2}+\dfrac{1}{2}\Vert\bar{\gamma}-\mathcal{P}_{0}\bar{\gamma}\Vert_{0,\Omega}^{2}+\dfrac{\alpha^{2}}{2}\Vert\gamma_{0}-\mathcal{P}_{0}(\gamma_{0})\Vert_{0,\Omega}^{2}\\
\delta\Vert e_{1}\Vert_{0,\Omega}^{2}  &  \leq Ch^{2}\left(  \Vert\bfu_{\varepsilon}\cdot\bfv_{\varepsilon}\Vert_{1,\Omega}^{2}+\Vert\bar{\gamma}\Vert_{1,\Omega}^{2}+\alpha^{2}\Vert\gamma_{0}\Vert_{1,\Omega}^{2}\right),
\end{align*}
where, applying H\"{o}lder inequality, Lemma \ref{adjoint}, Theorems \ref{f-bound} and \ref{H2}, we have
\begin{align*}
\Vert\bfu_{\varepsilon}\cdot\bfv_{\varepsilon}\Vert_{1,\Omega}  & =\Vert(\nabla\bfu_{\varepsilon})\bfv_{\varepsilon}+(\nabla\bfv_{\varepsilon})\bfu_{\varepsilon}\Vert_{0,\Omega}\\
&  \leq\Vert(\nabla\bfu_{\varepsilon})\bfv_{\varepsilon}\Vert_{0,\Omega}+\Vert(\nabla\bfv_{\varepsilon})\bfu_{\varepsilon}\Vert_{0,\Omega}\\
&  \leq\Vert\bfu_{\varepsilon}\Vert_{1,\Omega}\Vert\bfv_{\varepsilon}\Vert_{0,\infty,\Omega}+\Vert\bfv_{\varepsilon}\Vert_{1,\Omega}\Vert\bfu_{\varepsilon}\Vert_{0,\infty,\Omega}\\
&  \leq C(\Vert\bfu_{0}\Vert_{0,\omega}+\Vert\bff\Vert_{0,\Omega})^{2},
\end{align*}
and
\begin{align*}
\delta\Vert e_{1}\Vert_{0,\Omega}^{2}  &  \leq Ch^{2}\left(  \Vert\bfu_{\varepsilon}^{h}\cdot\bfv_{\varepsilon}^{h}\Vert_{1,\Omega}^{2}+\Vert\bar{\gamma}\Vert_{1,\Omega}^{2}+\alpha^{2}\Vert\gamma_{0}\Vert_{1,\Omega}^{2}\right) \\
\Vert e_{1}\Vert_{0,\Omega}  &  \leq\dfrac{C}{\delta^{1/2}}h((\Vert\bfu_{0}\Vert_{0,\omega}+\Vert\bff\Vert_{0,\Omega})^{2}+\Vert\bar{\gamma}\Vert_{1,\Omega}+\alpha\Vert\gamma_{0}\Vert_{1,\Omega})\\
&  \leq\dfrac{C}{\delta^{1/2}}h((\Vert\bfu_{0}\Vert_{0,\omega}+\Vert\bff\Vert_{0,\Omega})^{2}+\alpha\Vert\gamma_{0}\Vert_{1,\Omega}).
\end{align*}
Let $e_{2}=\gamma_{\varepsilon}^{h}-\hat{\gamma}_{\varepsilon}^{h}$. By Taylor Theorem, there exists $t\in[0,1]$ and $\xi_{t}=t\gamma_{\varepsilon}^{h}+(1-t)\hat{\gamma}_{\varepsilon}^{h}\in\mathcal{A}_{\varepsilon}^{h}(\bar{\gamma})$ such that%
\[
J_{h}^{\prime\prime}(\xi_{t})\left[  e_{2},e_{2}\right]  =J_{h}^{\prime}(\gamma_{\varepsilon}^{h})[e_{2}]-J_{h}^{\prime}(\hat{\gamma}_{\varepsilon}^{h})[e_{2}],
\]
with
\[
J_{h}^{\prime}(\gamma_{\varepsilon}^{h})[e_{2}]=J_{h}^{\prime}(\gamma_{\varepsilon}^{h})[\gamma_{\varepsilon}^{h}-\hat{\gamma}_{\varepsilon}^{h}]\leq0\leq J^{\prime}(\hat{\gamma}_{\varepsilon}^{h})[\gamma_{\varepsilon}^{h}-\hat{\gamma}_{\varepsilon}^{h}]=J^{\prime}(\hat{\gamma}_{\varepsilon}^{h})[e_{2}].
\]
Then,%
\begin{align*}
\dfrac{\delta}{4}\Vert e_{2}\Vert_{0,\Omega}^{2}  &  \leq J_{h}^{\prime\prime}(\xi_{t})\left[  e_{2},e_{2}\right]  =J_{h}^{\prime}(\gamma_{\varepsilon}^{h})[e_{2}]-J_{h}^{\prime}(\hat{\gamma}_{\varepsilon}^{h})[e_{2}]\\
&  \leq J^{\prime}(\hat{\gamma}_{\varepsilon}^{h})[e_{2}]-J_{h}^{\prime}(\hat{\gamma}_{\varepsilon}^{h})[e_{2}]\\
&  \leq Ch^{2}\Vert e_{2}\Vert_{0,\Omega}(\Vert\bff\Vert_{0,\Omega}+\Vert\bfu_{0}\Vert_{0,\omega})^{2}\\
\Vert e_{2}\Vert_{0,\Omega}  &  \leq\dfrac{C}{\delta}h^{2}(\Vert\bff\Vert_{0,\Omega}+\Vert\bfu_{0}\Vert_{0,\omega})^{2}.
\end{align*}
In conclusion,%
\[
\Vert\bar{\gamma}-\gamma_{\varepsilon}^{h}\Vert_{0,\Omega}\leq\Vert e_{1}\Vert_{0,\Omega}+\Vert e_{2}\Vert_{0,\Omega}\leq\dfrac{C}{\delta^{1/2}}h((\Vert\bfu_{0}\Vert_{0,\omega}+\Vert\bff\Vert_{0,\Omega})^{2}+\alpha\Vert\gamma_{0}\Vert_{1,\Omega})
\]
proving this result.
\end{proof}

Finally, we present our a priori error estimate for this first case.
\begin{theorem}\label{conv-p0}
Let $\bar{\gamma}\in\mathcal{A}$ a local solution of \eqref{min} that satisfies the optimality conditions \eqref{first-order} and \eqref{second-order}. For $h_0>0$ small enough, there exists a sequence $\{\gamma^{h}\}_{0<h<h_0}\ $ of solutions to the discrete optimal control problem \eqref{min-h} and a positive constant $C$ such that for all $h\in (0,h_0)$
\[
\Vert\bar{\gamma}-\gamma^{h}\Vert_{0,\Omega}\leq\dfrac{C}{\delta^{1/2}}h((\Vert\bfu_{0}\Vert_{0,\omega}+\Vert\bff\Vert_{0,\Omega})^{2}+\alpha\Vert\gamma_{0}\Vert_{1,\Omega}).
\]
\end{theorem}

\begin{proof}
First, there exists $\varepsilon>0$ small enough such that the solution $\bar{\gamma}_{\varepsilon}^{h}$ of the auxiliary optimal control problem \eqref{aux-oc} is also a local solution for the discrete optimal control problem \eqref{min-h} for $h>0$ small enough. The estimate is a consequence of Theorem \ref{P0}.
\end{proof}

\subsection{Case with a continuous Lagrange discrete control}

In this subsection, we consider the case $\mathcal{A}^{h}=\mathcal{A}\cap G_{1}^{h}$, also called $\PP_1$ case. The deduction of our a priori error estimates is slightly different from the previous case. Let $\bar{\gamma}\in\mathcal{A}$ be a local solution of \eqref{min}, the set $\mathcal{T}_{h}$ can be partitioned as $\mathcal{T}_{h}=\mathcal{T}_{h}^{1}\cup\mathcal{T}_{h}^{2}\cup\mathcal{T}_{h}^{3}$, where
\begin{align*}
\mathcal{T}_{h}^{1}  &  =\left\{  T\in\mathcal{T}_{h}\text{ }\mid\text{ }\bar{\gamma}\vert_{T}=a\text{ or }\bar{\gamma}\vert_{T}=b\right\}, \\
\mathcal{T}_{h}^{2}  &  =\left\{  T\in\mathcal{T}_{h}\text{ }\mid\text{}a<\bar{\gamma}<b\text{ on }T\right\}, \\
\mathcal{T}_{h}^{3}  &  =\mathcal{T}_{h}^{1}\setminus(\mathcal{T}_{h}^{2}\cup\mathcal{T}_{h}^{3}).
\end{align*}
We denote $\Omega_{h,j}=\bigcup\limits_{T\in\mathcal{T}_{h}^{j}}T$ for $j\in\left\{  1,2,3\right\}  $ and we assume that $\vert\Omega_{h,3}\vert\leq Ch^{p}$ for some $p\geq1$. This kind of hypothesis has been discussed in \cite{KV09,MR04}.

\begin{theorem}
Consider $\bar{\gamma}$ a local solution of \eqref{min} that satisfies the optimality conditions \eqref{first-order} and \eqref{second-order}. For $h_0>0$ small enough, there exists a sequence $\{\gamma^{h}\}_{0<h<h_0}\ $ of solutions to the discrete optimal control problem \eqref{min-h} and a positive constant $C$ such that for all $h\in (0,h_0)$
\[
\Vert\bar{\gamma}-\gamma^{h}\Vert_{0,\Omega}\leq\left(  1+\dfrac{C}{\delta}\right)  \Vert\bar{\gamma}-\mathcal{I}_{1}(\bar{\gamma})\Vert_{0,\Omega}+\dfrac{C}{\delta^{1/2}}\sqrt{J^{\prime}(\bar{\gamma})[\mathcal{I}_{1}(\bar{\gamma})-\bar{\gamma}]}+\dfrac{C}{\delta}h^{2}.
\]
\end{theorem}

\begin{proof}
As in the proof of Theorem \ref{P0}, there exists a constant $\varepsilon>0$ such that for all $\varphi\in L^{\infty}(\Omega)$ and all $\gamma\in\mathcal{A}_{\varepsilon}(\bar{\gamma})$
\[
J^{\prime\prime}(\gamma)[\varphi,\varphi]\geq\dfrac{\delta}{2}\Vert\varphi\Vert_{0,\Omega}^{2}
\]
and a constant $h_0>0$ such that for all $h\in(0,h_0)$, $\varphi\in L^{\infty}(\Omega)$ and $\gamma\in\mathcal{A}_{\varepsilon}(\bar{\gamma})$%
\[
J_{h}^{\prime\prime}(\gamma)[\varphi,\varphi]\geq\dfrac{\delta}{4}\Vert\varphi\Vert_{0,\Omega}^{2}.
\]
Then, there exists an unique solution $\hat{\gamma}_{\varepsilon}^{h} \in\mathcal{A}_{\varepsilon}^{h}(\bar{\gamma})$ for the auxiliary optimal control problem \eqref{aux-oc}. Now, we have%
\[
\Vert\bar{\gamma}-\gamma^{h}\Vert_{0,\Omega}\leq\Vert\bar{\gamma}-\mathcal{I}_{1}(\bar{\gamma})\Vert_{0,\Omega}+\Vert\mathcal{I}_{1}(\bar{\gamma})-\gamma^{h}\Vert_{0,\Omega}.
\]
Let $e=\bar{\gamma}-\gamma^{h}$, $e_{1}=\bar{\gamma}-\mathcal{I}_{1}(\bar{\gamma})$ and $e_{2}=\mathcal{I}_{1}(\bar{\gamma})-\gamma^{h}$. Then,
\[
J^{\prime}(\bar{\gamma})[e]=J^{\prime}(\bar{\gamma})[\bar{\gamma}-\gamma^{h}]\leq0\leq J_{h}^{\prime}(\gamma^{h})[\mathcal{I}_{1}(\bar{\gamma})-\gamma^{h}]=J_{h}^{\prime}(\gamma^{h})[e_{2}].
\]
By Taylor Theorem, there exists $t\in[0,1]$ and $\xi_{t}=t\mathcal{I}_{1}(\bar{\gamma})+(1-t)\gamma^{h}\in\mathcal{A}_{\varepsilon}(\bar{\gamma})$ such that
\begin{align*}
J_{h}^{\prime\prime}(\xi_{t})\left[  e_{2},e_{2}\right]  =  &  J_{h}^{\prime}(\mathcal{I}_{1}(\bar{\gamma}))[e_{2}]-J_{h}^{\prime}(\gamma^{h})[e_{2}]\\
\leq &  J_{h}^{\prime}(\mathcal{I}_{1}(\bar{\gamma}))[e_{2}]-J^{\prime}(\bar{\gamma})[e]\\
\leq &  (J_{h}^{\prime}(\mathcal{I}_{1}(\bar{\gamma}))[e_{2}]-J_{h}^{\prime}(\bar{\gamma})[e_{2}])+(J_{h}^{\prime}(\bar{\gamma})[e_{2}]-J^{\prime}(\bar{\gamma})[e_{2}])-J^{\prime}(\bar{\gamma})[e_{1}],
\end{align*}
where, applying Lemma \ref{aux-1} and Proposition \ref{interpolator},
\begin{align*}
J_{h}^{\prime}(\mathcal{I}_{1}(\bar{\gamma}))[e_{2}]-J_{h}^{\prime}(\bar{\gamma})[e_{2}]  &  \leq C\Vert e_{2}\Vert_{0,\Omega}\Vert\mathcal{I}_{1}(\bar{\gamma})-\bar{\gamma}\Vert_{0,\Omega}=C\Vert e_{2}\Vert_{0,\Omega}\Vert e_{1}\Vert_{0,\Omega}\\
J_{h}^{\prime}(\bar{\gamma})[e_{2}]-J^{\prime}(\bar{\gamma})[e_{2}]  &  \leq Ch^{2}\Vert e_{2}\Vert_{0,\Omega}.
\end{align*}
Then,%
\begin{align*}
\dfrac{\delta}{4}\Vert e_{2}\Vert_{0,\Omega}^{2}  &  \leq J_{h}^{\prime\prime}(\xi_{t})\left[  e_{2},e_{2}\right]  \leq C\Vert e_{2}\Vert_{0,\Omega}(h^{2}+\Vert e_{1}\Vert_{0,\Omega})-J^{\prime}(\bar{\gamma})[e_{1}]\\
&  \leq C\Vert e_{2}\Vert_{0,\Omega}(h^{2}+\Vert e_{1}\Vert_{0,\Omega})+J^{\prime}(\bar{\gamma})[\mathcal{I}_{1}(\bar{\gamma})-\bar{\gamma}],
\end{align*}
proving that%
\[
\Vert e_{2}\Vert_{0,\Omega}\leq\dfrac{C}{\delta}(h^{2}+\Vert e_{1}\Vert_{0,\Omega})+\dfrac{C}{\delta}\sqrt{J^{\prime}(\bar{\gamma})[\mathcal{I}_{1}(\bar{\gamma})-\bar{\gamma}]}.
\]
Finally,%
\begin{align*}
\Vert\bar{\gamma}-\gamma^{h}\Vert_{0,\Omega}  &  \leq\Vert\bar{\gamma}-\mathcal{I}_{1}(\bar{\gamma})\Vert_{0,\Omega}+\Vert\mathcal{I}_{1}(\bar{\gamma})-\gamma^{h}\Vert_{0,\Omega}\\
&  \leq\Vert\bar{\gamma}-\mathcal{I}_{1}(\bar{\gamma})\Vert_{0,\Omega}+\dfrac{C}{\delta}(h^{2}+\Vert\bar{\gamma}-\mathcal{I}_{1}(\bar{\gamma})\Vert_{0,\Omega})+\dfrac{C}{\delta}\sqrt{J^{\prime}(\bar{\gamma})[\mathcal{I}_{1}(\bar{\gamma})-\bar{\gamma}]}\\
&  \leq C\left(  1+\dfrac{1}{\delta}\right)  \Vert\bar{\gamma}-\mathcal{I}_{1}(\bar{\gamma})\Vert_{0,\Omega}+\dfrac{C}{\delta^{1/2}}\sqrt{J^{\prime}(\bar{\gamma})[\mathcal{I}_{1}(\bar{\gamma})-\bar{\gamma}]}+\dfrac{C}{\delta}h^{2},
\end{align*}
proving the theorem.
\end{proof}

The assumption $\vert\Omega_{h,3}\vert\leq Ch^{p}$ is necessary to obtain a new interpolation estimate.
\begin{lemma}\label{aux-p1}
Suppose $\gamma_{0}\in H^{2}(\Omega)\cap W^{1,\infty}(\Omega)$ and $\bar{\gamma}\in W^{1,\infty}(\Omega)$. Then, $\bar{\gamma}\in H^{2}(T)$ for each $T\in\mathcal{T}_{h}^{2}$. Furthermore, there exists a constant $C>0$ such that
\[
\Vert\bar{\gamma}-\mathcal{I}_{1}(\bar{\gamma})\Vert_{0,\Omega}^{2}=C(h^{4}\Vert\triangle\bar{\gamma}\Vert_{0,\Omega_{h,2}}^{2}+h^{2+p}\Vert\nabla\bar{\gamma}\Vert_{0,\infty,\Omega}^{2}).
\]

\end{lemma}

\begin{proof}
First, if $A(\bar{\gamma})=(\boldsymbol{\bar{u}},\bar{p})$ and $B(\bar{\gamma})=(\boldsymbol{\bar{v}},\bar{q})$, we have $\bar{\gamma}=\gamma_{0}+\dfrac{1}{\alpha}(\boldsymbol{\bar{u}}\cdot\boldsymbol{\bar{v}})$ for all $T\in\mathcal{T}_{h}^{2}$, proving that $\bar{\gamma}\in H^{2}(T)$. Second, we have
\[
\Vert\bar{\gamma}-\mathcal{I}_{1}(\bar{\gamma})\Vert_{0,\Omega}^{2}=\sum\limits_{T\in\mathcal{T}_{h}^{1}}\Vert\bar{\gamma}-\mathcal{I}_{1}(\bar{\gamma})\Vert_{0,T}^{2}+\sum\limits_{T\in\mathcal{T}_{h}^{2}}\Vert\bar{\gamma}-\mathcal{I}_{1}(\bar{\gamma})\Vert_{0,T}^{2}+\sum\limits_{T\in\mathcal{T}_{h}^{3}}\Vert\bar{\gamma}-\mathcal{I}_{1}(\bar{\gamma})\Vert_{0,T}^{2},
\]
where $\bar{\gamma}=\mathcal{I}_{1}(\bar{\gamma})$ on $T\in\mathcal{T}_{h}^{1}$. Then, applying the local Lagrange interpolator estimates (Theorem 1.103 in \cite{EG04}),
\[
\sum\limits_{T\in\mathcal{T}_{h}^{2}}\Vert\bar{\gamma}-\mathcal{I}_{1}(\bar{\gamma})\Vert_{0,T}^{2}\leq C\sum\limits_{T\in\mathcal{T}_{h}^{2}}h_{T}^{4}\Vert\triangle\bar{\gamma}\Vert_{0,T}^{2}\leq Ch^{4}\Vert\triangle\bar{\gamma}\Vert_{0,\Omega_{h,2}}^{2}.
\]
If $T\in\mathcal{T}_{h}^{3}$, we have that $a=\min\limits_{x\in K}\bar{\gamma}(x)<\max\limits_{x\in K}\bar{\gamma}(x)<b$ or $a<\min\limits_{x\in K}\bar{\gamma}(x)<\max\limits_{x\in K}\bar{\gamma}(x)=b$. In the first case, we have%
\[
\Vert\bar{\gamma}-\mathcal{I}_{1}(\bar{\gamma})\Vert_{0,T}^{2}\leq\Vert T\Vert\Vert\bar{\gamma}-\mathcal{I}_{1}(\bar{\gamma})\Vert_{0,\infty,T}^{2}\leq\Vert T\Vert\Vert\bar{\gamma}-a\Vert_{0,\infty,T}^{2}\leq\Vert T\Vert h_{T}^{2}\Vert\nabla\bar{\gamma}\Vert_{0,\infty,T}^{2},
\]
as a consequence of the Mean Value Theorem. We obtain a similar estimate for the second case. In consequence,%
\[
\sum\limits_{T\in\mathcal{T}_{h}^{3}}\Vert\bar{\gamma}-\mathcal{I}_{1}(\bar{\gamma})\Vert_{0,T}^{2}\leq\sum\limits_{T\in\mathcal{T}_{h}^{3}}\Vert T\Vert h_{T}^{2}\Vert\nabla\bar{\gamma}\Vert_{0,\infty,T}^{2}\leq Ch^{2+p}\Vert\nabla\bar{\gamma}\Vert_{0,\infty,\Omega_{h,3}}^{2},
\]
proving the desired estimate.
\end{proof}

Our a priori error estimate for this case is given in the following theorem.
\begin{theorem}\label{conv-p1}
Consider $\bar{\gamma}\in W^{1,\infty}(\Omega)$ a local solution of \eqref{min} that fulfills the optimality conditions \eqref{first-order} and \eqref{second-order}. For $h_0>0$ small enough, there exist a positive constant $C$ and a sequence $\{\gamma^{h}\}_{0<h<h_0}$ of solutions of the discrete optimal control problem \eqref{min-h} such that for all $h\in (0,h_0)$
\[
\Vert\bar{\gamma}-\gamma^{h}\Vert_{0,\Omega}\leq\left(  1+\dfrac{C}{\delta}\right)  h^{1+p/2}+\dfrac{C}{\delta}h^{2}.
\]
\end{theorem}

\begin{proof}
Consider $A(\bar{\gamma})=(\boldsymbol{\bar{u}},\bar{p})$ and $B(\bar{\gamma})=(\boldsymbol{\bar{v}},\bar{q})$. Since $\mathcal{I}_{1}(\bar{\gamma})=\bar{\gamma}$ in $\Omega_{h,1}$ and $\bar{\gamma}=\gamma_{0}+\dfrac{1}{\alpha}(\boldsymbol{\bar{u}}\cdot\boldsymbol{\bar{v}})$, we have%
\[
0\leq J^{\prime}(\bar{\gamma})[\mathcal{I}_{1}(\bar{\gamma})-\bar{\gamma}]=\sum\limits_{T\in\mathcal{T}_{h}^{3}}(\mathcal{I}_{1}(\bar{\gamma})-\bar{\gamma},\boldsymbol{\bar{u}}\cdot\boldsymbol{\bar{v}}+\alpha(\bar{\gamma}-\gamma_{0}))_{T}.
\]
Defining $d=\boldsymbol{\bar{u}}\cdot\boldsymbol{\bar{v}}+\alpha(\bar{\gamma}-\gamma_{0})$, we have that for all $T\in\mathcal{T}_{h}^{3}$, there exists $\boldsymbol{x}_{T}\in T$ such that $d(\boldsymbol{x}_{T})=\boldsymbol{0}$. Repeating the same reasoning as in the proof of Lemma \ref{aux-p1}, we obtain
\begin{align*}
(\mathcal{I}_{1}(\bar{\gamma})-\bar{\gamma},\boldsymbol{\bar{u}}\cdot\boldsymbol{\bar{v}}+\alpha(\bar{\gamma}-\gamma_{0}))_{T}  & =(\mathcal{I}_{1}(\bar{\gamma})-\bar{\gamma},d-d(\boldsymbol{x}_{T}))_{T}\\
&  \leq\Vert\mathcal{I}_{1}(\bar{\gamma})-\bar{\gamma}\Vert_{0,T}\Vert d-d(\boldsymbol{x}_{T})\Vert_{0,T}\\
&  \leq\vert T\vert h_{T}^{2}\Vert\nabla\bar{\gamma}\Vert_{0,\infty,T}^{2}.
\end{align*}
Then,%
\[
J^{\prime}(\bar{\gamma})[\mathcal{I}_{1}(\bar{\gamma})-\bar{\gamma}]=\sum\limits_{T\in\mathcal{T}_{h}^{3}}(\mathcal{I}_{1}(\bar{\gamma})-\bar{\gamma},\boldsymbol{\bar{u}}\cdot\boldsymbol{\bar{v}}+\alpha(\bar{\gamma}-\gamma_{0}))_{T}\leq\sum\limits_{T\in\mathcal{T}_{h}^{3}}\Vert T\Vert h_{T}^{2}\Vert\nabla\bar{\gamma}\Vert_{0,\infty,T}^{2}\leq Ch^{2+p}\Vert\nabla\bar{\gamma}\Vert_{0,\infty,\Omega_{h,3}}^{2},
\]
proving that%
\[
\Vert\bar{\gamma}-\gamma^{h}\Vert_{0,\Omega}\leq\left(  1+\dfrac{C}{\delta}\right)  h^{1+p/2}+\dfrac{C}{\delta}h^{2}.
\]
\end{proof}

\subsection{Semi-discrete scheme}
The optimal control problem can also be defined with a discretization of the states and respective adjoint, but with no discretizations for the control $\gamma$. In this case, the local solutions $\bar{\gamma}^{h}\in\mathcal{A}$ must verify the identity $\bar{\gamma}^{h}=\Pi_{[ a,b]}\left(\gamma_{0}+\dfrac{1}{\alpha}\boldsymbol{\bar{u}}^{h}\cdot\boldsymbol{\bar{v}}^{h}\right)  $, where $(\boldsymbol{\bar{u}}^{h},\bar{p}^{h})$ and $(\boldsymbol{\bar{v}}^{h},\bar{q}^{h})$ are the respective optimal states and adjoints. The a priori error estimates are almost a direct consequence of the results obtained in Section 2.

\begin{theorem}\label{semi}
Let $\bar{\gamma}\in\mathcal{A}$ a local solution of \eqref{min} that verifies the optimality conditions \eqref{first-order} and \eqref{second-order}. For $h_{0}>0$ small enough, there exists a sequence $\{\gamma^{h}\}_{0<h<h_{0}}\in\mathcal{A}$ of solutions of the semi-discrete optimal control problem \eqref{semi-h} such that
\[
\Vert\gamma^{h}-\bar{\gamma}\Vert_{0,\Omega}\leq C\Vert\bff%
\Vert_{0,\Omega}(\Vert\bff\Vert_{0,\Omega}+\Vert\bfu%
_{0}\Vert_{0,\omega}).
\]

\end{theorem}

\begin{proof}
From Lemma \ref{ellipticity}, there exists $\varepsilon>0$ such that for all $\varphi\in L^{\infty}(\Omega)$ and all $\gamma\in\mathcal{A}_{\varepsilon}(\bar{\gamma})$
\[
J^{\prime\prime}(\gamma)[\varphi,\varphi]\geq\dfrac{\delta}{2}\Vert\varphi\Vert_{0,\Omega}^{2}.
\]
Then, the semi-discrete optimal control problem has an unique solution $\gamma^{h}\in\mathcal{A}_{\varepsilon}(\bar{\gamma})$ for $h>0$ small enough. Taking $\varphi=\gamma^{h}-\bar{\gamma}$ and applying Taylor theorem, there exists $t\in[0,1]$ and $\xi_{t}=t\gamma^{h}+(1-t)\bar{\gamma}\in\mathcal{A}_{\varepsilon}(\bar{\gamma})$ such that
\[
\dfrac{\delta}{2}\Vert\gamma^{h}-\bar{\gamma}\Vert_{0,\Omega}^{2}\leq J^{\prime\prime}(\xi_{t})[\gamma^{h}-\bar{\gamma},\gamma^{h}-\bar{\gamma}]=J^{\prime}(\gamma^{h})[\gamma^{h}-\bar{\gamma}]-J(\bar{\gamma})[\gamma^{h}-\bar{\gamma}].
\]
Since $\bar{\gamma}$ and $\gamma^{h}$ fulfill their respective optimality conditions, we have
\[
J^{\prime}(\bar{\gamma})\left[  \bar{\gamma}-\hat{\gamma}_{\varepsilon}^{h}\right]  \leq0\leq J_{h}^{\prime}(\gamma^{h})[\gamma^{h}-\bar{\gamma}]
\]
Then, applying Lemma \ref{aux-1},
\begin{align*}
\dfrac{\delta}{2}\Vert\gamma^{h}-\bar{\gamma}\Vert_{0,\Omega}^{2}  &  \leq J^{\prime\prime}(\xi_{t})[\gamma^{h}-\bar{\gamma},\gamma^{h}-\bar{\gamma}]=J^{\prime}(\gamma^{h})[\gamma^{h}-\bar{\gamma}]-J(\bar{\gamma})[\gamma^{h}-\bar{\gamma}]\\
&  \leq J^{\prime}(\gamma^{h})[\gamma^{h}-\bar{\gamma}]-J_{h}^{\prime}(\gamma^{h})[\gamma^{h}-\bar{\gamma}]\\
&  \leq Ch^{2}\Vert\gamma^{h}-\bar{\gamma}\Vert_{0,\Omega}\Vert\bff\Vert_{0,\Omega}(\Vert\bff\Vert_{0,\Omega}+\Vert\bfu_{0}\Vert_{0,\omega}),
\end{align*}
proving that $\Vert\gamma^{h}-\bar{\gamma}\Vert_{0,\Omega}\leq Ch^{2}\Vert\bff\Vert_{0,\Omega}(\Vert\bff\Vert_{0,\Omega}+\Vert\bfu_{0}\Vert_{0,\omega})$.
\end{proof}

\begin{remark}
The results obtained in Theorems \ref{taylor-hood}, \ref{conv-p0}, \ref{conv-p1} and \ref{semi} can be extended to nonhomogeneous Dirichlet boundary conditions given by $\bfu = \bfg$ on $\partial\Omega$ if $\bfg\in\bfH^{1/2}(\partial\Omega)$ satisfies the compatibility condition $\langle \bfg,\bfn \rangle_{1/2,\partial\Omega} =0$. In that case, the Dirichlet boundary condition can be analyzed following same arguments that we use in, for example, a Poisson problem with a nonhomogeneous Dirichlet boundary condition.
\end{remark}

\section{A posteriori error estimate}\label{S4}
In this section we discuss the deduction of an a posteriori error estimator. First we introduce our estimators. 
\begin{definition}
Let $\bar{\gamma}\in\mathcal{A}$ and $\gamma^{h}\in\mathcal{A}$ local solutions of \eqref{min} and \eqref{min-h}, respectively, such that $A(\bar{\gamma})=(\boldsymbol{\bar{u}},\bar{p})$, $B(\bar{\gamma})=(\boldsymbol{\bar{v}},\bar{q})$, $A_{h}(\gamma^{h})=(\bfu^{h},p^{h})$ and $B(\gamma^{h})=(\bfv^{h},q^{h})$. We define $e(\gamma)=\bar{\gamma}-\gamma^{h}$, $e(\bfu)=\boldsymbol{\bar{u}}-\bfu^{h}$, $e(p)=\bar{p}-p^{h}$,$e(\bfv)=\boldsymbol{\bar{v}}-\bfv^{h}\ $and $e(q)=\bar{q}-q^{h}$. For $T\in\mathcal{T}_{h}$ and $F\in\mathcal{E}_{h}$, we denote
\begin{align*}
\mathcal{R}_{T}  &  =\bff+\nu\triangle\bfu^{h}-(\nabla\bfu^{h})\bfu^{h}-\nabla p^{h}-\gamma^{h}\bfu^{h},\\
\mathcal{R}_{A,T}  &  =\chi_{\omega}(\bfu^{h}-\bfu_{0})+\nu\triangle\bfv^{h}+(\nabla\bfv^{h})\bfu^{h}-(\nabla\bfu^{h})^{T}\bfv^{h}-\nabla q^{h}-\gamma^{h}\bfv^{h},\\
\mathcal{J}_{F}  &  =\jump{(\nu\nabla\bfu^{h}-p^{h}\boldsymbol{I})\boldsymbol{n}},\\
\mathcal{J}_{A,F}  &  =\jump{(\nu\nabla\bfv^{h}+q^{h}\boldsymbol{I})\boldsymbol{n}},
\end{align*}
where $\chi_\omega\in L^{\infty}(\Omega)$ is the indicator function for $\omega$ and $\jump{\cdot}$ denotes the jump operator (see Equation 1.54 in \cite{EG04}). For $\gamma^{\ast}=\Pi_{[ a,b]}\left(\gamma_{0}+\dfrac{1}{\alpha}(\boldsymbol{\bar{u}}^{h}\cdot\boldsymbol{\bar{v}}^{h})\right)$, we define $(\bfu^{\ast},p^{\ast})=A(\gamma^{\ast})$ and $(\bfv^{\ast},q^{\ast})=B(\gamma^{\ast})$. Finally, we define
\begin{align*}
\eta_{S,T}^{2}  &  =h_{T}^{2}\Vert\mathcal{R}_{T}\Vert_{0,T}^{2}+h_{T}\Vert\mathcal{J}_{F}\Vert_{0,\partial T\setminus\partial\Omega}^{2}+\Vert\operatorname{div}\bfu^{h}\Vert_{0,T}^{2},\\
\eta_{A,T}^{2}  &  =h_{T}^{2}\Vert\mathcal{R}_{A,T}\Vert_{0,T}^{2}+h_{T}\Vert\mathcal{J}_{A,F}\Vert_{0,\partial T\setminus\partial\Omega}^{2}+\Vert\operatorname{div}\bfv^{h}\Vert_{0,T}^{2},\\
\eta_{C,T}^{2}  &  =\Vert\gamma^{h}-\gamma^{\ast}\Vert_{0,T}^{2},\\
\eta_{S,h}^{2}  &  =\sum_{T\in\mathcal{T}_{h}}\eta_{S,T}^{2}\quad \eta_{A,h}^{2}=\sum_{T\in\mathcal{T}_{h}}\eta_{A,T}^{2}\quad \eta
_{C,h}^{2}=\sum_{T\in\mathcal{T}_{h}}\eta_{C,T}^{2},\\
\eta_{T}^{2}  &  =\left\{\begin{array}
[c]{ll}%
\eta_{S,T}^{2}+\eta_{A,T}^{2}+\eta_{C,T}^{2} & \text{for the discrete scheme}\\
\eta_{S,T}^{2}+\eta_{A,T}^{2} & \text{for the semi-discrete scheme}
\end{array}
\right. \\
\eta_{h}^{2}  &  =\sum_{T\in\mathcal{T}_{h}}\eta_{T}^{2}=\left\{
\begin{array}
[c]{ll}%
\eta_{S,h}^{2}+\eta_{A,h}^{2}+\eta_{C,h}^{2} & \text{for the discrete
scheme}\\
\eta_{S,h}^{2}+\eta_{A,h}^{2} & \text{for the semi-discrete scheme}%
\end{array}
\right.
\end{align*}
\end{definition}
The estimators $\eta_S$ and $\eta_A$ are the same residual a posteriori error estimators for the Navier-Stokes and Oseen equations, respectively, defined in \cite{OWA94} and \cite{V96}. Now, we present an auxiliary result.
\begin{lemma}\label{aux-rel}
Under the same hypotheses of Theorem \ref{semi}, for all $\gamma\in\mathcal{A}$, $(\alpha(\gamma^{\ast}-\gamma_{0})+\bfu^{h}\cdot\bfv^{h},\bar{\gamma}-\gamma^{\ast})_\Omega\geq 0$.
\end{lemma}
\begin{proof}
See Lemma 2.26 in \cite{T10}.
\end{proof}
\subsection{Reliability of the a posterior error estimators}
In this subsection, we present an upper bound for each component of the total a posteriori error estimator, considering the control, the states and the adjoints errors. First, we present two classical results from \cite{OWA94} and \cite{V96}.
\begin{lemma}\label{rel-states}
There exists a constant $C>0$ independent on $h$ such that $\Vert A(\gamma^{h})-A_{h}(\gamma^{h})\Vert\leq C\eta_{S,h}$ and $\Vert B(\gamma^{h})-B_{h}(\gamma^{h})\Vert\leq C\eta_{A,h}$.
\end{lemma}
\begin{proof}
See \cite{OWA94} and \cite{V96}.
\end{proof}
Now we can prove the reliability of our a posterior error estimators that is valid for our three schemes.
\begin{theorem}
There exists a constant $C>0$ independent on $h$ such that
\begin{equation*}
\Vert(e(\bfu),e(p))\Vert+\Vert(e(\bfv),e(q))\Vert+\Vert e(\gamma)\Vert_{0,\Omega}\leq C\eta_{h}.
\end{equation*}
\end{theorem}

\begin{proof}
First, we present the proof for the discrete scheme. We have%
\[
\Vert e(\gamma)\Vert_{0,\Omega}=\Vert\bar{\gamma}-\gamma^{h}\Vert_{0,\Omega}\leq\Vert\bar{\gamma}-\gamma^{\ast}\Vert_{0,\Omega}+\Vert\gamma^{\ast}-\gamma^{h}\Vert_{0,\Omega}=\Vert\bar{\gamma}-\gamma^{\ast}\Vert_{0,\Omega}+\eta_{C,h},
\]
and we define $(\boldsymbol{\hat{u}},\hat{p})=A(\gamma^{h})$ and $(\boldsymbol{\hat{v}},\hat{q})=B(\gamma^{h})$. Reasoning as in the proof of Lemma \ref{taylor}, taking $\varphi=\bar{\gamma}-\gamma^{\ast}$ and applying Taylor theorem, there exists $t\in [0,1]$ and $\xi_{t}=t\bar{\gamma}+(1-t)\gamma^{\ast}\in\mathcal{A}_{\varepsilon}(\bar{\gamma})$ such that%
\[
\dfrac{\delta}{2}\Vert\bar{\gamma}-\gamma^{\ast}\Vert_{0,\Omega}^{2}\leq J^{\prime\prime}(\xi_{t})[\bar{\gamma}-\gamma^{\ast},\bar{\gamma}-\gamma^{\ast}]=J^{\prime}(\bar{\gamma})[\bar{\gamma}-\gamma^{\ast}]-J^{\prime}(\gamma^{\ast})[\bar{\gamma}-\gamma^{\ast}],
\]
where $J^{\prime}(\bar{\gamma})[\bar{\gamma}-\gamma^{\ast}]\leq 0$. Applying Lemma \ref{aux-rel},
\begin{align*}
\dfrac{\delta}{2}\Vert\gamma^{h}-\bar{\gamma}\Vert_{0,\Omega}^{2}  & \leq(\alpha(\gamma^{\ast}-\gamma_{0})+\bfu^{h}\cdot\bfv^{h},\bar{\gamma}-\gamma^{\ast})_\Omega-J^{\prime}(\gamma^{\ast})[\bar{\gamma}-\gamma^{\ast}]\\
&  \leq(\bfu^{h}\cdot\bfv^{h}-\bfu^{\ast}\cdot\bfv^{\ast},\bar{\gamma}-\gamma^{\ast})_\Omega\\
&  \leq\Vert\gamma^{h}-\bar{\gamma}\Vert_{0,\Omega}\Vert\bfu^{h}\cdot\bfv^{h}-\bfu^{\ast}\cdot\bfv^{\ast}\Vert_{0,\Omega}.
\end{align*}
Then, by H\"older inequality and Sobolev Embedding Theorem,%
\begin{align*}
\Vert\gamma^{h}-\bar{\gamma}\Vert_{0,\Omega}  &  \leq C\Vert\bfu^{h}\cdot\bfv^{h}-\bfu^{\ast}\cdot\bfv^{\ast}\Vert_{0,\Omega}\\
&  \leq C(\Vert\bfu^{h}\Vert_{0,\Omega}\Vert\bfv^{h}-\bfv^{\ast}\Vert_{0,\Omega}+\Vert\bfu^{h}-\bfu^{\ast}\Vert_{0,\Omega}\Vert\bfv^{\ast}\Vert_{0,\Omega})\\
&  \leq C(\Vert\bfu^{h}\Vert_{0,4,\Omega}\Vert\bfv^{h}-\bfv^{\ast}\Vert_{0,4,\Omega}+\Vert\bfu^{h}-\bfu^{\ast}\Vert_{0,4,\Omega}\Vert\bfv^{\ast}\Vert_{0,4,\Omega})\\
&  \leq C(\vert\bfu^{h}\vert_{1,\Omega}\vert\bfv^{h}-\bfv^{\ast}\vert_{1,\Omega}+\vert\bfu^{h}-\bfu^{\ast}\vert_{1,\Omega}\vert\boldsymbol{\hat{v}}\vert_{1,\Omega}),
\end{align*}
where, applying Lemmas \ref{lipschitz} and \ref{rel-states},%
\begin{align*}
\vert\bfu^{h}-\bfu^{\ast}\vert_{1,\Omega}  &  \leq\vert \bfu^{h}-\boldsymbol{\hat{u}}\vert_{1,\Omega}+\vert \boldsymbol{\hat{u}}-\bfu^{\ast}\vert_{1,\Omega}\\
&  \leq\Vert A(\gamma^{h})-A_{h}(\gamma^{h})\Vert+\Vert A(\gamma^{h})-A(\gamma^{\ast})\Vert\\
&  \leq C\eta_{S,h}+C\Vert\gamma^{h}-\gamma^{\ast}\Vert=C(\eta_{S,h}+\eta_{C,h}).
\end{align*}
Analogously, $\vert \bfu^{h}-\bfu^{\ast}\vert_{1,\Omega}\leq C(\eta_{A,h}+\eta_{C,h})$. Since $\vert \bfu^{h}\vert_{1,\Omega}\leq C\Vert\bff\Vert_{0,\Omega}$ and $\vert \boldsymbol{\hat{v}}\vert_{1,\Omega}\leq C(\Vert\bff\Vert_{0,\Omega}+\Vert\bfu_{0}\Vert_{0,\omega})$, we obtain
\[
\Vert\gamma^{h}-\bar{\gamma}\Vert_{0,\Omega}\leq C(\eta_{S,h}+\eta_{A,h}+\eta_{C,h})=C\eta_{h}%
\]
proving that%
\[
\Vert e(\gamma)\Vert_{0,\Omega}\leq\Vert\bar{\gamma}-\gamma^{\ast}\Vert_{0,\Omega}+\eta_{C,T}\leq C\eta_{h}%
\]
Following the same reasoning, we have
\begin{align*}
\Vert(e(\bfu),e(p))\Vert &  \leq\Vert A(\bar{\gamma})-A(\gamma^{h})\Vert+\Vert A(\gamma^{h})-A_{h}(\gamma^{h})\Vert\\
&  \leq C\Vert\bar{\gamma}-\gamma^{h}\Vert_{0,\Omega}+C\eta_{S,h}\\
&  \leq C\eta_{h},
\end{align*}
and%
\begin{align*}
\Vert(e(\bfv),e(q))\Vert &  \leq\Vert B(\bar{\gamma})-B(\gamma^{h})\Vert+\Vert B(\gamma^{h})-B_{h}(\gamma^{h})\Vert\\
&  \leq C\Vert\bar{\gamma}-\gamma^{h}\Vert_{0,\Omega}+C\eta_{A,h}\\
&  \leq C\eta_{h}.
\end{align*}
In conclusion, $\Vert(e(\bfu),e(p))\Vert+\Vert(e(\bfv),e(q))\Vert+\Vert e(\gamma)\Vert_{0,\Omega}\leq C\eta_{h}$. For the semi-discrete scheme, we have $\gamma^{h}=\gamma^{\ast}$. Then, we repeat the same previous analysis to obtain the estimate. We omit the details for this case.
\end{proof}

\subsection{Efficiency of the a posteriori error estimators}
In this subsection, we present a lower bound for each term of our local a posteriori error estimators. First, the $\eta_C$ term is analyzed. 
\begin{lemma}\label{eff-control}
There exists a constant $C>0$ independent on $h$ such that
\begin{equation*}
\eta_{C,h}\leq \Vert e(\gamma)\Vert_{0,\Omega} + C(\Vert(e(\bfu),e(p))\Vert+\Vert(e(\bfv),e(q))\Vert).
\end{equation*}
\end{lemma}

\begin{proof}
We have%
\[
\eta_{C,h}=\Vert\gamma^{\ast}-\gamma_{h}\Vert_{0,\Omega}\leq\Vert\gamma^{\ast}-\bar{\gamma}\Vert_{0,\Omega}+\Vert\bar{\gamma}-\gamma^{h}\Vert_{0,\Omega}=\Vert\gamma^{\ast}-\bar{\gamma}\Vert_{0,\Omega}+\Vert e(\gamma)\Vert_{0,\Omega}.
\]
Since $\Pi_{[ a,b]}(\cdot)$ is Lipschitz, applying H\"older inequality and Sobolev Embedding Theorem,
\begin{align*}
\Vert\gamma^{\ast}-\bar{\gamma}\Vert_{0,\Omega}  &  =\left\Vert \Pi_{[a,b]}\left(  \gamma_{0}+\dfrac{1}{\alpha}\bfu^{h}\cdot\bfv^{h}\right)  -\Pi_{[ a,b]}\left(  \gamma_{0}+\dfrac{1}{\alpha}\boldsymbol{\bar{u}}\cdot\boldsymbol{\bar{v}}\right)  \right\Vert_{0,\Omega}\\
&  \leq\dfrac{1}{\alpha}\Vert\bfu^{h}\cdot\bfv^{h}-\boldsymbol{\bar{u}}\cdot\boldsymbol{\bar{v}}\Vert_{0,\Omega}\\
&  \leq\dfrac{1}{\alpha}\Vert\bfu^{h}\cdot(\bfv^{h}-\boldsymbol{\bar{v}})+(\bfu^{h}-\boldsymbol{\bar{u}})\cdot\boldsymbol{\bar{v}}\Vert_{0,\Omega}\\
&  \leq\dfrac{1}{\alpha}(\Vert\bfu^{h}\Vert_{0,4,\Omega}\Vert\bfv^{h}-\boldsymbol{\bar{v}}\Vert_{0,4,\Omega}+\Vert\bfu^{h}-\boldsymbol{\bar{u}}\Vert_{0,4,\Omega}\Vert\boldsymbol{\bar{v}}\Vert_{0,4,\Omega})\\
&  \leq\dfrac{C}{\alpha}(\vert \bfu^{h}\vert_{1,\Omega}\vert \bfv^{h}-\boldsymbol{\bar{v}}\vert_{1,\Omega}+\vert \bfu^{h}-\boldsymbol{\bar{u}}\vert_{1,\Omega}\vert \boldsymbol{\bar{v}}\vert_{1,\Omega})\\
&  \leq\dfrac{C}{\alpha}(\Vert\bff\Vert_{0,\Omega}\Vert(e(\bfu),e(p))\Vert+(\Vert\bff\Vert_{0,\Omega}+\Vert\bfu_{0}\Vert_{0,\omega})\Vert(e(\bfv),e(q))\Vert)\\
&  \leq C(\Vert(e(\bfu),e(p))\Vert+\Vert(e(\bfv),e(q))\Vert),
\end{align*}
proving this lemma.
\end{proof}

Following the same reasoning as in \cite{V96}, we define the element and edge bubble functions and present some important properties.
\begin{definition}
For $T\in\mathcal{T}_{h}$ and $E\in\mathcal{E}_{h}$, we denote $\omega(T)=\bigcup\{T^{\prime}\in\mathcal{T}_{h}$ $\mid$ $\overline{T}\cap\overline{T^{\prime}}\neq\emptyset\}$ and $\omega(E)=\bigcup\{T^{\prime}\in\mathcal{T}_{h}$ $\mid$ $E\cap\overline{T^{\prime}}\neq\emptyset\}$. We denote by $\psi_{T}$ and $\psi_{E}$ the element and edge bubble functions and by $\mathcal{P}:C(E)  \rightarrow C(T)  $ the continuation operator (see Section 3.1 in \cite{V96}).
\end{definition}

\begin{lemma}\label{bubble}
Given $k\in\NN$, there exists a constant $C>0$, depending only on $k$ and the shape-regularity of $\mathcal{T}_{h}$, such that for all $T\in\mathcal{T}_h$, all $E$ edge of $T$, all $q\in\PP_{k}\left(  T\right)  $, and all $r\in\PP_{k}\left(  E\right)  $ we have
\begin{align*}
\Vert q\Vert _{0,T}^{2}  &  \leq C\Vert \psi_{T}^{1/2}q\Vert _{0,T}^{2},\\
\Vert r\Vert _{0,E}^{2}  &  \leq C\Vert \psi_{E}^{1/2}r\Vert _{0,E}^{2},\\
\Vert \psi_{E}^{1/2} \mathcal{P}(r)\Vert _{0,T}^{2}  &  \leq Ch_{E}\Vert r\Vert _{0,E}^{2}.
\end{align*}
\end{lemma}
\begin{proof}
See Lemma 3.3 in \cite{V96}.
\end{proof}

Let $(\bfw,r)\in V$. Then, applying integration by parts on each $T\in\mathcal{T}_{h}$, we have%
\begin{align*}
&  a(e(\bfu),\bfw)+c(e(\bfu),\boldsymbol{\bar{u}},\bfw)+c(\bfu^{h},e(\bfu),\bfw)-b(\bfw,e(p))+b(e(\bfu),r)+(e(\gamma)\boldsymbol{\bar{u}},\bfw)+(\gamma^{h}e(\bfu),\bfw)_\Omega\\
= & \sum\limits_{T\in\mathcal{T}_{h}} [(\mathcal{R}_{T},\bfw)_{T}-(r,\operatorname{div}(\bfu^{h}))_\Omega]  +\sum_{E\in\mathcal{E}_{h}}(\mathcal{J}_{E},\bfw)_{E}
\end{align*}
and%
\begin{align*}
&  a(e(\bfv),\bfw)-b(\boldsymbol{r},e(q))+b(e(\bfv),r)+c(\boldsymbol{\bar{u}},\bfw,e(\bfv))+c(e(\bfu),\bfw,\bfv^{h})+(e(\gamma)\boldsymbol{\bar{v}},\bfw)_\Omega+(\gamma^{h}e(\boldsymbol{\bar{v}}),\bfw)_\Omega\\
= &\sum\limits_{T\in\mathcal{T}_{h}} [(\mathcal{R}_{A,T},\bfw)_{T}-(r,\operatorname{div}(\bfv^{h}))_\Omega] +\sum_{E\in\mathcal{E}_{h}}(\mathcal{J}_{A,E},\bfw)_{E}%
\end{align*}
Now we proceed by cases. For the discrete scheme, we define
\begin{align*}
\mathcal{R}_{T}^{0}  &  =P(\bff)+\nu\triangle\bfu^{h}-(\nabla\bfu^{h})\bfu^{h}-\nabla p^{h}-\gamma^{h}\bfu^{h},\\
\mathcal{R}_{A,T}^{0}  &  =\chi_{\omega}\bfu^{h}-P(\chi_{\omega}\bfu_{0})+\nu\triangle\bfv^{h}+(\nabla\bfv^{h})\bfu^{h}-(\nabla\bfu^{h})^{T}\bfv^{h}-\nabla q^{h}-\gamma^{h}\bfv^{h}.
\end{align*}
Then,
\begin{align*}
\sum\limits_{T\in\mathcal{T}_{h}}(\mathcal{R}_{T},\bfw)_T  & =\sum\limits_{T\in\mathcal{T}_{h}}(\mathcal{R}_{T}^{0},\bfw)_T+\sum\limits_{T\in\mathcal{T}_{h}}(\bff-P(\bff),\bfw)_T,\\
\sum\limits_{T\in\mathcal{T}_{h}}(\mathcal{R}_{A,T},\bfw)_T  & =\sum\limits_{T\in\mathcal{T}_{h}}(\mathcal{R}_{T}^{0},\bfw)_T+\sum\limits_{T\in\mathcal{T}_{h}}(\chi_{\omega}\bfu_{0}-P(\chi_{\omega}\bfu_{0}),\bfw)_T.
\end{align*}
Similarly, we define for the semi-discrete scheme%
\begin{align*}
\mathcal{R}_{T}^{0}  &  =P(\bff)+\nu\triangle\bfu^{h}-(\nabla\bfu^{h})\bfu^{h}-\nabla p^{h}-P(\bar{\gamma}\boldsymbol{\bar{u}}),\\
\mathcal{R}_{A,T}^{0}  &  =\chi_{\omega}\bfu^{h}-P(\chi_{\omega}\bfu_{0})+\nu\triangle\bfv^{h}+(\nabla\bfv^{h})\bfu^{h}-(\nabla\bfu^{h})^{T}\bfv^{h}-\nabla q^{h}-P(\bar{\gamma}\boldsymbol{\bar{v}}),
\end{align*}
where%
\begin{align*}
\sum\limits_{T\in\mathcal{T}_{h}}(\mathcal{R}_{T},\bfw)_T  & =\sum\limits_{T\in\mathcal{T}_{h}}\left[  (\mathcal{R}_{T}^{0},\bfw)_T+(e(\gamma)\bfu^{h},\bfw)_T+(\bar{\gamma}e(\bfu),\bfw)_T+(\bff-P(\bff),\bfw)_T-(\bar{\gamma}\boldsymbol{\bar{u}}-P(\bar{\gamma}\boldsymbol{\bar{u}}),\bfw)_T\right], \\
\sum\limits_{T\in\mathcal{T}_{h}}(\mathcal{R}_{A,T},\bfw)_T  & =\sum\limits_{T\in\mathcal{T}_{h}}\left[  (\mathcal{R}_{T}^{0},\bfw)_T+(e(\gamma)\bfv^{h},\bfw)_T+(\bar{\gamma}e(\bfv),\bfw)_T+(\chi_{\omega}\bfu_{0}-P(\chi_{\omega}\bfu_{0}),\bfw)_T-(\bar{\gamma}\boldsymbol{\bar{v}}-P(\bar{\gamma}\boldsymbol{\bar{v}}),\bfw)_T\right].
\end{align*}

\begin{definition}
Let $\boldsymbol{g}\in\boldsymbol{L}^{2}(\Omega)$ we denote by $\Theta(\bff)$ the oscillation residual term given by%
\[
\Theta(\boldsymbol{g})=\left(  \sum\limits_{T\in\mathcal{T}_{h}}h_{T}^{2}%
\Vert\boldsymbol{g}-P(\boldsymbol{g})\Vert_{0,T}^{2}\right)  ^{1/2}.
\]
\end{definition}

First, we detail the proof of the efficiency for the discrete cases.

\begin{lemma}\label{eff-states}
For the discrete cases, there exists a constant $C>0$ independent on $h$ such that
\begin{equation*}
\eta_{S,h}\leq C\left(  \Vert(e(\bfu),e(p))\Vert+\Vert e(\gamma)\Vert_{0,\Omega}+\Theta(\bff)\right).
\end{equation*}
\end{lemma}

\begin{proof}
First, consider $T\in\mathcal{T}_{h}$. Taking $(\bfw,r)=(\psi_{T}\mathcal{R}_{T}^{0},0)$, since $\bfw=\boldsymbol{0}$ in $\Omega\setminus T$,
\begin{align*}
\Vert\psi_{T}^{1/2}\mathcal{R}_{T}^{0}\Vert_{0,T}^{2}=  &  \sum\limits_{T\in\mathcal{T}_{h}}(\mathcal{R}_{T}^{0},\bfw)_T-\sum\limits_{T\in\mathcal{T}_{h}}(\bff-P(\bff),\bfw)_T\\
=  &  a(e(\bfu),\bfw)+c(e(\bfu),\boldsymbol{\bar{u}},\bfw)+c(\bfu^{h},e(\bfu),\bfw)-b(\bfw,e(p))\\
&  +(e(\gamma)\boldsymbol{\bar{u}},\bfw)_\Omega+(\gamma^{h}e(\bfu),\bfw)_\Omega-\sum\limits_{T\in\mathcal{T}_{h}}(\bff-P(\bff),\bfw)_T.
\end{align*}
Applying Lemmas \ref{bounds-1} and \ref{bubble} , H\"older inequality and Sobolev Embedding theorem,%
\begin{align*}
\Vert\psi_{T}^{1/2}\mathcal{R}_{T}^{0}\Vert_{0,T}^{2}\leq &  \nu\vert e(\bfu)\vert_{1,T}\vert \bfw\vert_{1,T}+\beta\vert \bfw\vert_{1,T} \vert e(\bfu)\vert_{1,T}(\vert \boldsymbol{\bar{u}}\vert_{1,T}+\vert \bfu^{h}\vert_{1,T})+\sqrt{d}\vert \boldsymbol{w \vert }_{1,T}\Vert e(p)\Vert_{0,T}\\
&  +\Vert e(\gamma)\Vert_{0,T}\Vert\boldsymbol{\bar{u}}\cdot\bfw\Vert_{0,T}+\Vert\gamma^{h}\Vert_{0,T}\Vert e(\bfu)\cdot\bfw\Vert_{0,T}+\Vert\bff-P(\bff)\Vert_{0,T}\Vert\bfw\Vert_{0,T}\\
\leq &  \nu \vert e(\bfu)\vert_{1,T}\vert \bfw\vert_{1,T} +\beta\vert \bfw\vert_{1,T} \vert e(\bfu)\vert_{1,T}(\vert \boldsymbol{\bar{u}}\vert_{1,T}+\vert \bfu^{h}\vert_{1,T})+\sqrt{d}\vert \boldsymbol{w \vert }_{1,T}\Vert e(p)\Vert_{0,T}\\
&  +\Vert e(\gamma)\Vert_{0,T}\Vert\boldsymbol{\bar{u}}\Vert_{0,4,T} \Vert\bfw\Vert_{0,T}+\Vert\gamma^{h}\Vert_{0,T}\Vert e(\bfu)\Vert_{0,4,T}\Vert\bfw\Vert_{0,4,T}+\Vert\bff-P(\bff)\Vert_{0,T}\Vert\bfw\Vert_{0,T}\\
\leq &  \nu \vert e(\bfu)\vert_{1,T}\vert \bfw\vert_{1,T}+\beta\vert \bfw\vert_{1,T} \vert e(\bfu)\vert_{1,T}(\vert \boldsymbol{\bar{u}}\vert_{1,T}+\vert \bfu^{h}\vert_{1,T})+\sqrt{d}\vert \bfw\vert_{1,T}\Vert e(p)\Vert_{0,T}\\
&  +C\vert \bfw\vert_{1,T}(\Vert e(\gamma)\Vert_{0,T}\vert \boldsymbol{\bar{u}}\vert_{1,T}+\Vert\gamma^{h}\Vert_{0,T} \vert e(\bfu)\vert_{1,T})+\Vert\bff-P(\bff)\Vert_{0,T}\Vert\bfw\Vert_{0,T}\\
\leq &  C\vert \bfw\vert_{1,T}( \vert e(\bfu)\vert_{1,T}+\Vert e(p)\Vert_{0,T}+\Vert e(\gamma)\Vert_{0,T})+\Vert\bff-P(\bff)\Vert_{0,T}\Vert\bfw\Vert_{0,T}.
\end{align*}
Since $\bfw$ is polynomial, we have $\vert \bfw\vert_{1,T}\leq Ch_{T}^{-1}\Vert\bfw\Vert_{0,T}\leq Ch_{T}^{-1}\Vert\mathcal{R}_{T}^{0}\Vert_{0,T}$ by applying an inverse inequality (see Lemma 1.138 from \cite{EG04}). Then,
\begin{align*}
C\Vert\mathcal{R}_{T}^{0}\Vert_{0,T}^{2}  &  \leq\Vert\psi_{T}^{1/2}\mathcal{R}_{T}^{0}\Vert_{0,T}^{2}\leq C\vert \bfw\vert_{1,T}(\vert e(\bfu)\vert_{1,T}+\Vert e(p)\Vert_{0,T}+\Vert e(\gamma)\Vert_{0,T})+\Vert\bff-P(\bff)\Vert_{0,T}\Vert\bfw\Vert_{0,T}\\
\Vert\mathcal{R}_{T}^{0}\Vert_{0,T}  &  \leq Ch_{T}^{-1}( \vert e(\bfu)\vert_{1,T}+\Vert e(p)\Vert_{0,T}+\Vert e(\gamma)\Vert_{0,T})+C\Vert\bff-P(\bff)\Vert_{0,T}\\
h_{T}^{2}\Vert\mathcal{R}_{T}^{0}\Vert_{0,T}^{2}  &  \leq C( \vert e(\bfu)\vert_{1,T}^{2}+\Vert e(p)\Vert_{0,T}^{2}+\Vert e(\gamma)\Vert_{0,T}^{2})+Ch_{T}^{2}\Vert\bff-P(\bff)\Vert_{0,T}^{2},
\end{align*}
but $\Vert\mathcal{R}_{T}\Vert_{0,T}\leq\Vert\mathcal{R}_{T}^{0}\Vert _{0,T}+\Vert\bff-P(\bff)\Vert_{0,T}$. Then,
\begin{align}
h_{T}^{2}\Vert\mathcal{R}_{T}\Vert_{0,T}^{2}  &  \leq C( \vert e(\bfu)\vert_{1,T}^{2}+\Vert e(p)\Vert_{0,T}^{2}+\Vert e(\gamma)\Vert_{0,T}^{2}+h_{T}^{2}\Vert\bff-P(\bff)\Vert_{0,T}^{2})\label{eff-1}\\
\sum\limits_{T\in\mathcal{T}_{h}}h_{T}^{2}\Vert\mathcal{R}_{T}\Vert_{0,T}^{2} &  \leq C\left(  \Vert e(\bfu),e(p)\Vert+\Vert e(\gamma)\Vert_{0,\Omega}^{2}+\sum\limits_{T\in\mathcal{T}_{h}}h_{T}^{2}\Vert\bff-P(\bff)\Vert_{0,T}^{2}\right).\nonumber
\end{align}
Now consider $E\in\mathcal{E}_{h}$ and $(\bfw,r)=(\psi_{E}P_{E}(\mathcal{J}_{E}),0)$. Since $\bfw=\boldsymbol{0}$ in $\Omega\setminus\omega(E)$, reasoning as before, we have
\begin{align*}
\Vert\psi_{E}^{1/2}\mathcal{J}_{E}\Vert_{0,E}^{2}=  &  \sum_{E\in\mathcal{E}_{h}}(\mathcal{J}_{E},\bfw)_E\\
=  &  a(e(\bfu),\bfw)+c(e(\bfu),\boldsymbol{\bar{u}},\bfw)+c(\bfu^{h},e(\bfu),\bfw)-b(\bfw,e(p))\\
&  +(e(\gamma)\boldsymbol{\bar{u}},\bfw)_\Omega+(\gamma^{h}e(\bfu),\bfw)_\Omega-\sum\limits_{T\in\mathcal{T}_{h}}\left[  (\bff-P(\bff),\bfw)_T+(\mathcal{R}_{T},\bfw)_T\right] \\
\leq &  C\vert \bfw\vert_{1,\omega(E)}( \vert e(\bfu)\vert_{1,\omega(E)}+\Vert e(p)\Vert_{0,\omega(E)}+\Vert e(\gamma)\Vert_{0,\omega(E)})\\
& +\left( \sum\limits_{E\in\mathcal{E}_{h}}\Vert\bff-P(\bff)\Vert_{0,\omega(E)}+\Vert\mathcal{R}_{T}\Vert_{0,\omega(E)}\right) \Vert\bfw\Vert_{0,\omega(E)},
\end{align*}
where $\vert \bfw\vert_{1,\omega(E)}\leq Ch_{E}^{-1}\Vert\bfw \Vert_{0,\omega(E)}\leq Ch_{E}^{-1/2}\Vert\mathcal{J}_{E}\Vert_{0,E}$ and $\Vert\mathcal{J}_{E}\Vert_{0,E}\leq C\Vert\psi_{E}^{1/2}\mathcal{J}_{E} \Vert_{0,E}$ by applying some inverse inequalities (see Lemma 1.138 in \cite{EG04}). Then,
\begin{align}
C\Vert\mathcal{J}_{E}\Vert_{0,E}^{2}  &  \leq\Vert\psi_{E}^{1/2}\mathcal{J}_{E}\Vert_{0,E}^{2}\nonumber\\
& \leq C \vert \bfw\vert_{1,\omega(E)}( \vert e(\bfu)\vert_{1,\omega(E)}+\Vert e(p)\Vert_{0,\omega(E)}+\Vert e(\gamma)\Vert_{0,\omega(E)})+(\Vert\bff-P(\bff)\Vert_{0,\omega(E)}+\Vert\mathcal{R}_{T}\Vert_{0,\omega(E)})\Vert
\bfw\Vert_{0,\omega(E)}\nonumber\\
h_{E}^{1/2}\Vert\mathcal{J}_{E}\Vert_{0,E}  &  \leq C( \vert e(\bfu)\vert_{1,\omega(E)}+\Vert e(p)\Vert_{0,\omega(E)}+\Vert e(\gamma)\Vert_{0,\omega(E)}+h_{E}(\Vert\bff-P(\bff)\Vert_{0,\omega(E)}+\Vert\mathcal{R}_{T}\Vert_{0,\omega(E)}))\nonumber\\
h_{E}\Vert\mathcal{J}_{E}\Vert_{0,E}^{2}  &  \leq C( \vert e(\bfu)\vert_{1,\omega(E)}^{2}+\Vert e(p)\Vert_{0,\omega(E)}^{2}+\Vert e(\gamma)\Vert_{0,\omega(E)}^{2}+h_{E}^{2}(\Vert\bff-P(\bff)\Vert_{0,\omega(E)}^{2}+\Vert\mathcal{R}_{T}\Vert_{0,\omega(E)}^{2}))\label{eff-2}\\
\sum_{E\in\mathcal{E}_{h}}h_{E}\Vert\mathcal{J}_{E}\Vert_{0,E}^{2}  &  \leq C\sum_{E\in\mathcal{E}_{h}}( \vert e(\bfu)\vert_{1,\omega(E)}^{2}+\Vert e(p)\Vert_{0,\omega(E)}^{2}+\Vert e(\gamma)\Vert_{0,\omega(E)}^{2}+h_{E}^{2}(\Vert\bff-P(\bff)\Vert_{0,\omega(E)}^{2}+\Vert\mathcal{R}_{T}\Vert_{0,\omega(E)}^{2}))\nonumber\\
\sum_{E\in\mathcal{E}_{h}}h_{E}\Vert\mathcal{J}_{E}\Vert_{0,E}^{2}  &  \leq C\left(  \Vert e(\bfu),e(p)\Vert^{2}+\Vert e(\gamma)\Vert_{0,\Omega}^{2}+\sum\limits_{T\in\mathcal{T}_{h}}h_{T}^{2}\Vert\bff-P(\bff)\Vert_{0,T}^{2}\right).\nonumber
\end{align}
Finally, it is direct that
\begin{equation}
\Vert\operatorname{div}\bfu^{h}\Vert_{0,T}^{2}=\Vert\operatorname{div}(\boldsymbol{\bar{u}}-\bfu^{h})\Vert_{0,T}^{2}\leq\sqrt{d} \vert \boldsymbol{\bar{u}}-\bfu^{h}\vert_{1,T}^{2}. \label{eff-3}
\end{equation}
Then, $\sum\limits_{T\in\mathcal{T}_{h}}\Vert\operatorname{div}\bfu^{h}\Vert_{0,T}^{2}\leq\sqrt{d}\Vert(e(\bfu),e(p))\Vert^{2}$. The final estimate is obtained by minor algebraic manipulations.
\end{proof}
 
The estimates \eqref{eff-1}, \eqref{eff-2} and \eqref{eff-3} constitute a local efficiency property for estimator $\eta_{S,h}$. The same arguments can be applied to the estimator $\eta_{A,h}$. In what follows, we assume that $\overline{\omega}$ can be decomposes using elements from $\mathcal{T}_h$, i.e., there exists $\mathcal{T}_{\omega}\subseteq\mathcal{T}_{h}$ such that $\overline{\omega}=\bigcup\left\{  T^{\prime}\in\mathcal{T}_{\omega}\right\}$.
\begin{lemma}\label{eff-adj}
For both discrete cases, there exists a constant $C>0$ independent on $h$ such that
\begin{equation*}
\eta_{A,h}\leq C(\Vert(e(\bfu),e(p))\Vert+\Vert(e(\bfv),e(q))\Vert+\Vert e(\gamma)\Vert_{0,\Omega}+\Theta(\bff)+\Theta(\chi_{\omega_{0}}\bfu_{0})).
\end{equation*}
\end{lemma}

\begin{proof}
The deduction follows the same scheme as the proof of the previous lemma, including local efficiency estimates similar to \eqref{eff-1}, \eqref{eff-2} and \eqref{eff-3}. We omit the details.
\end{proof}

\begin{theorem}
For both discrete cases, there exists a constant $C>0$ independent on $h$ such that
\[
\eta_{h}\leq C(\Vert(e(\bfu),e(p))\Vert+\Vert(e(\bfv),e(q))\Vert+\Vert e(\gamma)\Vert_{0,\Omega}+\Theta(\bff%
)+\Theta(\chi_{\omega_{0}}\bfu_{0})).
\]
\end{theorem}

\begin{proof}
It is direct consequence of Lemmas \ref{eff-control}, \ref{eff-states} and \ref{eff-adj}, and some algebraic manipulations. We omit the details.
\end{proof}

This theorem is only a global effectivity result for the estimator $\eta_h$ in both discrete cases. For the semi-discrete case, we can obtain similar estimates to \eqref{eff-1}, \eqref{eff-2} and \eqref{eff-3}, concluding a local efficiency property for estimator $\eta_{S,h}$ and $\eta_{S,h}$. Then, the estimator $\eta_h$ for the semi-discrete discretization is globally and locally efficient since $\eta_h$ does not include the term $\eta_{C,h}$. This result is summarized in the following theorem.

\begin{theorem}
For the semi-discrete scheme, there exists a constant $C>0$ independent on $h$ such that
\[
\eta_{h}\leq C(\Vert(e(\bfu),e(p))\Vert+\Vert(e(\bfv),e(q))\Vert+\Vert e(\gamma)\Vert_{0,\Omega}+\Theta(\bff)+\Theta(\chi_{\omega_{0}}\bfu_{0})+\Theta(\bar{\gamma}\boldsymbol{\bar{u}})+\Theta(\bar{\gamma}\boldsymbol{\bar{v}})).
\]
\end{theorem}

\begin{proof}
Following the same steps as in the proof of Lemmas \ref{eff-states} and \ref{eff-adj}, we can obtain the inequalities
\begin{align*}
\eta_{S,h}  &  \leq C\left(  \Vert(e(\bfu),e(p))\Vert+\Vert e(\gamma)\Vert_{0,\Omega}+\Theta(\bff)+\Theta(\bar{\gamma}\boldsymbol{\bar{u}})\right), \\
\eta_{A,h}  &  \leq C(\Vert(e(\bfu),e(p))\Vert+\Vert (e(\bfv),e(q))\Vert+\Vert e(\gamma)\Vert_{0,\Omega}+\Theta(\bff)+\Theta(\chi_{\omega_{0}}\bfu_{0})+\Theta (\bar{\gamma}\boldsymbol{\bar{u}})+\Theta(\bar{\gamma}\boldsymbol{\bar{v}})).
\end{align*}
We omit the details.
\end{proof}

\section{Numerical experiments}\label{S5}
In this section, we report some numerical experiments in 2D using the a priori and a posteriori error estimates established in this paper, and for each of the three schemes presented. For the case $\mathcal{A}^h=\mathcal{A}\cap G_{0}^{h}$, we solve the coupled variational formulation system given by \eqref{VFh}, \eqref{AVFh} and \eqref{first-order-p0}, replacing $\gamma$ with $\gamma^h$ for the first two variational formulations, using a non-linear semi-smooth Newton solver. For case $\mathcal{A}=\mathcal{A}\cap G_{1}^{h}$, we solve the same coupled variational formulation, but replacing the last formulation with the condition $\gamma^{h}=\mathcal{I}_{1}^{h}(\gamma^{\ast})$, using Picard method. The main reason is that the interpolation, using which the control can verify the box constraint, cannot be easily rewritten as a variational identity. Moreover, in this setting, the Picard method is restricted to large enough values of $\alpha$, similar to the real-world simulation.
Finally, the semi-discrete scheme can be rewritten as a coupled variational formulation given by \eqref{VFh} and \eqref{AVFh}, but replacing $\gamma$ with $\gamma^\ast$ and is solved using the same semi-smooth non-linear Newton solver as in the case $\mathcal{A}=G_{0}^{h}$. Conditions $(p,1)_\Omega=(q,1)_\Omega=0$ are imposed by two Lagrange multipliers.

Our examples were implemented in \texttt{FEniCS 2019.1.0} \cite{fenics} with a PETSc module. The semi-smooth non-linear Newton solver consists of a straightforward application of a Newton trust region algorithm \cite[Chapter 4]{NW06}, where the Jacobian is computed by automatic differentiation, with absolute and relative tolerances equal to $10^{-12}$. The tolerance for the Picard method is given by $10^{-6}$ for the discrete $l^2$ norm of two consecutive values of $\gamma^h$. Linear systems were solved using \texttt{MUMPS} \cite{MUMPS}, a direct solver suitable for large linear systems with sparse matrices. The numerical tests are run on an Intel Core i7-10750H @ 2.60 GHz computer, using Ubuntu 24.04.5 LTS as VM inside a Windows Subsystem for Linux (WSL2) with 32 GB of RAM.

We denote by \texttt{dof} the number of degrees of freedom of the discrete systems and we define a quotient called the effectivity index $\theta$ by
\[
\theta:=\dfrac{\eta}{(\Vert(e(\bfu),e(p))\Vert^{2}+\Vert(e(\bfv),e(q))\Vert ^{2}+\Vert e(\gamma)\Vert_{0,\Omega}^{2})^{1/2}},
\]
which, ideally, converges to $1$ as $h\rightarrow 0$. The numerical results have shown that this can be archived even on coarse meshes, indicating estimated errors to be close to true errors effectively. \cite{BR79} and \cite{BMV83} have shown that such an equilibrium of the error indicators will lead to an optimal mesh. However, this is not always observed in practice. As shown in \cite{ZZ87}, where their estimator is considered effective despite the fact that counterexamples exist where the effectivity index may not be asymptotically exact ($\theta$ converging to $1$). Similarly, \cite{LP05} has shown that in certain situations, for example, using low-order finite elements where aspect ratios are high, the effectivity index can be greater than 1 or even increase. This is observed, for example, in the semi-discrete test case, where $\eta^2_{C,h}$ and thus the penalization constant $\frac{1}{\alpha}$ vanishes. 

The algorithm of the implementation follows the standard AFEM procedure from most well-known literature. Originally introduced in \cite{D96} and further developed in \cite{CKNS08,MNS00}, a simple one-step AFEM is a loop of the sequence:
\[ 
\texttt{SOLVE} \longrightarrow \texttt{ESTIMATE} \longrightarrow \texttt{MARK} \longrightarrow \texttt{REFINE} 
\]

The solution will be computed on the current mesh in step \texttt{SOLVE}. Then using the local a posteriori error estimator $\eta_{T}$ for all $T\in\mathcal{T}_{h}$ to identify regions with large errors, our algorithm marks selected elements for refinement $T$ such that $\eta_{T}\geq\rho\max\{\eta_{\tilde{T}}$~$\mid$~$\tilde{T}$ $\in \mathcal{T}_{h}\}$ for a constant $\rho\in(0,1)$, this is also called Dörfler marking strategy \cite{D96}. Recent work in \cite{PP19} has developed a new algorithm for Dörfler marking using a minimal set $T$ at linear costs. The normal Dörfler marking is sufficient for our test cases. In regions where the solution is smooth, the mesh can be coarsened to reduce computational cost. Once refined or coarsened accordingly, the new mesh is ready for the next iteration. This loop will be repeated until the stopping criterion is reached. The algorithm in pseudo-code\footnote{The entire FeniCS code is freely accessible via a request to \href{mailto:jaguayo@dim.uchile.cl}{jaguayo@dim.uchile.cl}} is summarized below:
    \begin{algorithm}[htpb]
    \caption{AFEM algorithm for Navier-Stokes-Brinkman optimal control problem}\label{afem}
    \begin{algorithmic}[1]
    \State \textbf{(INITIALIZE)} Set the parameters $\nu$, $\bff$, $\omega$, $\alpha$, $\bfu_0$, $\gamma_0$, $\mathcal{T}_h$, $\rho$, the counter $k := 0$ and the tolerances \texttt{tol}, \texttt{maxit} and \texttt{maxdof}.
    \State \textbf{(SOLVE)} Define the finite element spaces for $\mathcal{T}_h$, solve the optimal control problem \eqref{min-h} on $\mathcal{T}_h$ and the non-linear system \eqref{VFh}-\eqref{AVFh}.
    \State \textbf{(ESTIMATE)} Compute the error estimator $\eta_T$ for each $T \in \mathcal{T}_h$.
    \State \textbf{(MARK)} Mark a subset $\mathcal{M}_k \subseteq \mathcal{T}_h$ such that $\mathcal{M}_k$ contains the elements $T$ such that $\eta_{T}\geq\rho\max\{\eta_{\tilde{T}}$~$\mid$~$\tilde{T}$ $\in \mathcal{T}_{h}\}$
    \State \textbf{(REFINE)} Refine each $T \in \mathcal{M}_k$ following the algorithm described in \cite{PC00}.
    \State If $\eta>\texttt{tol}$, $k<\texttt{maxit}$, and $\texttt{dof}<\texttt{maxdof}$ set $ k := k+1$ then go to Step 2.
    \end{algorithmic}
    \end{algorithm}

\subsection{Test with analytic solution in 2D}
\begin{figure}[H]
\centering
\includegraphics[width=0.9\textwidth]{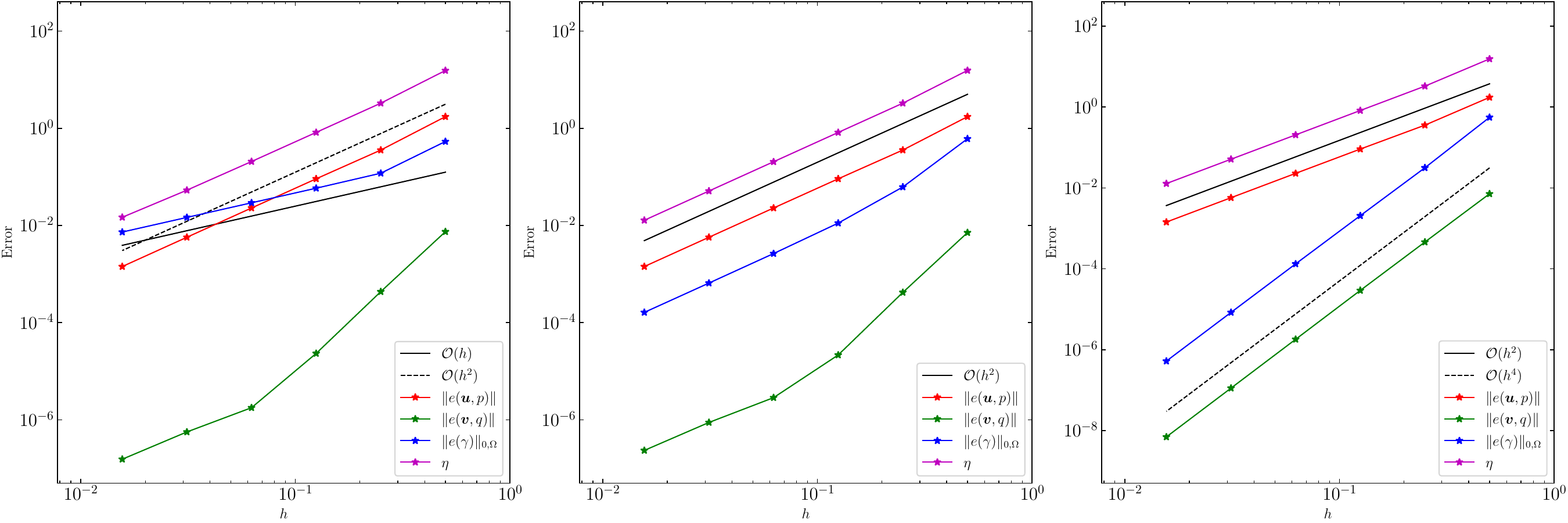}
\caption{History of convergence for $\mathcal{A}^h=\mathcal{A}\cap G_0 ^{h}$ (left), $\mathcal{A}^h=\mathcal{A}\cap G_1 ^{h}$ (center) and semi-discrete schemes (right).}
\label{error-ej1}
\end{figure}

In our first experiment, we consider the domain $\Omega=]-1,1[^{2}$ and the measurement subset $\omega=]-1/2,1/2[^{2}$. Taking $\bfu_{0}(x,y)=(\sin(\pi x)\sin(\pi y),\cos(\pi x)\cos(\pi y))^{T}$, $\gamma_{0}(x,y)=(1-x^{2})^{2}(1-y^{2})^{2}$, $a=0$, $b=1$, $\nu=1$ and $\alpha=10^{-3}$, we choose $\bff\in\bfL^{2}(\Omega)$ such that the unique solution of \eqref{min} is given by $\bar{\gamma}=\gamma_{0}$ and the optimal states and adjoints are given by
\begin{alignat*}{2}
\boldsymbol{\bar{u}}(x,y) &  =\bfu_{0}(x,y),  & \quad \bar{p}(x,y) & = xy,\\
\boldsymbol{\bar{v}}(x,y) &  =\boldsymbol{0}, & \quad \bar{q}(x,y) & = 0.
\end{alignat*}
We discretize the domain $\Omega$ using a structured triangular mesh with uniform refinements and we summarize the number of iterations of each solver, and the a priori and a posteriori error estimations for our three discrete schemes in Figure \ref{error-ej1} and Tables \ref{ej1-p0}, \ref{ej1-p1} and \ref{ej1-semi}. The convergence rates for controls and states are similar to the theoretical rates obtained in Theorems \ref{taylor-hood}, \ref{conv-p0}, \ref{conv-p1} and \ref{semi}, with indications of superconvergence in the adjoint states of all the analyzed schemes and in the control obtained in the semi-discrete scheme.

\begin{table}[H]
\centering
\begin{tabular}[c]{|c|c|c|c|c|c|c|c|}\hline
$h$ & \texttt{dof} & It. & $\Vert e(\gamma)\Vert_{0,\Omega}$ & $\Vert(e(\bfu),e(p))\Vert$ & $\Vert(e(\bfv),e(q))\Vert$ & $\eta$& $\theta$\\\hline
$1/2$  & $408$    & $8$  & $5.34929E-01$ & $1.73430E+00$ & $7.46280E-03$ & $1.54548E+01$ & $8.5153$\\
$1/4$  & $1448$   & $18$ & $1.18998E-01$ & $3.56553E-01$ & $4.35664E-04$ & $3.28459E+00$ & $8.7383$\\
$1/8$  & $5448$   & $21$ & $5.84945E-02$ & $9.12815E-02$ & $2.33276E-05$ & $8.22919E-01$ & $7.5904$\\
$1/16$ & $21128$  & $22$ & $2.93046E-02$ & $2.29105E-02$ & $1.78356E-06$ & $2.07039E-01$ & $5.5659$\\
$1/32$ & $83208$  & $23$ & $1.46603E-02$ & $5.72617E-03$ & $5.65451E-07$ & $5.32815E-02$ & $3.3853$\\
$1/64$ & $330248$ & $24$ & $7.33111E-03$ & $1.43055E-03$ & $1.55426E-07$ & $1.47568E-02$ & $1.9756$\\\hline
\end{tabular}
\caption{A priori and a posteriori error estimates with effectivity indexes, $\PP_0$ case.}
\label{ej1-p0}
\end{table}

\begin{table}[H]
\centering
\begin{tabular}
[c]{|c|c|c|c|c|c|c|c|}\hline
$h$ & \texttt{dof} & It. & $\Vert e(\gamma)\Vert_{0,\Omega}$ & $\Vert(e(\bfu),e(p))\Vert$ & $\Vert(e(\bfv),e(q))\Vert$ & $\eta$& $\theta$\\\hline
$1/2$  & $401$    & $5$  & $6.08433E-01$ & $1.73451E+00$ & $7.14075E-03$ & $1.54968E+01$ & $8.4307$\\
$1/4$  & $1401$   & $19$ & $6.19725E-02$ & $3.56750E-01$ & $4.19961E-04$ & $3.28357E+00$ & $9.0683$\\
$1/8$  & $5225$   & $20$ & $1.12007E-02$ & $9.12990E-02$ & $2.16744E-05$ & $8.20869E-01$ & $8.9241$\\
$1/16$ & $20169$  & $21$ & $2.64598E-03$ & $2.29139E-02$ & $2.87230E-06$ & $2.04961E-01$ & $8.8858$\\
$1/32$ & $79241$  & $23$ & $6.53875E-04$ & $5.72695E-03$ & $8.86575E-07$ & $5.12266E-02$ & $8.8871$\\
$1/64$ & $314121$ & $24$ & $1.63028E-04$ & $1.43074E-03$ & $2.35782E-07$ & $1.28074E-02$ & $8.8940$\\\hline
\end{tabular}
\caption{A priori and a posteriori error estimates with effectivity indexes, $\PP_1$ case.}
\label{ej1-p1}
\end{table}

\begin{table}[H]
\centering
\begin{tabular}
[c]{|c|c|c|c|c|c|c|c|}\hline
$h$ & \texttt{dof} & It. & $\Vert e(\gamma)\Vert_{0,\Omega}$ & $\Vert(e(\bfu),e(p))\Vert$ & $\Vert(e(\bfv),e(q))\Vert$ & $\eta$& $\theta$\\\hline
$1/2$  & $376$    & $10$ & $5.55824E-01$ & $1.73416E+00$ & $7.14835E-03$ & $1.54999E+01$ & $8.5114$\\
$1/4$  & $1320$   & $18$ & $3.14857E-02$ & $3.56441E-01$ & $4.58696E-04$ & $3.28313E+00$ & $9.1751$\\
$1/8$  & $4936$   & $16$ & $2.05686E-03$ & $9.12648E-02$ & $2.91629E-05$ & $8.20791E-01$ & $8.9912$\\
$1/16$ & $19080$  & $14$ & $1.32230E-04$ & $2.29076E-02$ & $1.80318E-06$ & $2.04942E-01$ & $8.9463$\\
$1/32$ & $75016$  & $11$ & $8.33954E-06$ & $5.72549E-03$ & $1.11880E-07$ & $5.12220E-02$ & $8.9463$\\
$1/64$ & $297480$ & $8$  & $5.21622E-07$ & $1.43039E-03$ & $6.95613E-09$ & $1.28062E-02$ & $8.9529$\\\hline
\end{tabular}
\caption{A priori and a posteriori error estimates with effectivity indexes, semi-discrete case.}
\label{ej1-semi}
\end{table}

\subsection{Test with a L-shaped domain}
\begin{figure}[H]
\centering
\includegraphics[width=0.90\textwidth]{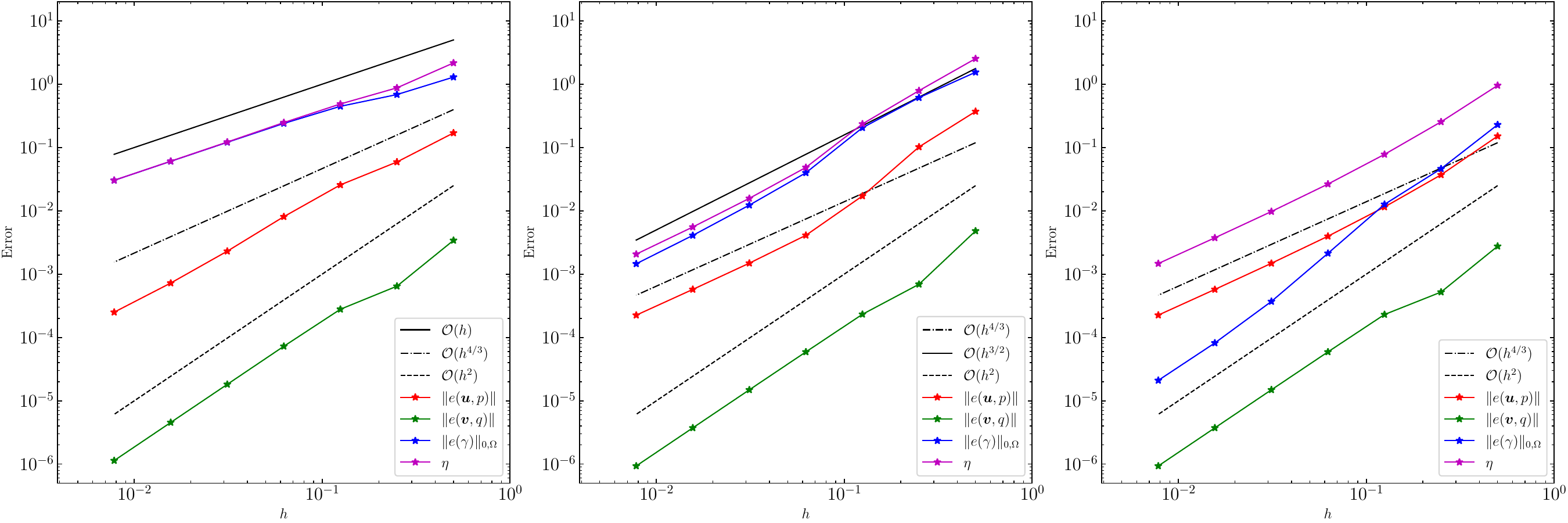}
\caption{History of convergence for $\PP_0$ (left), $\PP_1$ (center) and semi-discrete schemes (right).}
\label{error-ej2}
\end{figure}

We consider the domain $\Omega=]-1,1[^{2}\setminus ]-1,0[^2$ and the measurement region $\omega=\Omega$. Choosing $\gamma_0 =0$, $\nu=1$, $\alpha=10^{-4}$, $a=0$ and $b=5$, we set $\bfu_0$, $\bff$, and the Dirichlet boundary conditions such that the unique solution of \eqref{min} is given by
\begin{align*}
\bar{\bfu}(x,y)   & = (x+y)\exp\left(\dfrac{1}{2}(x+y)\right) (1,-1)^T,\\
\bar{p}(x,y)      & = r^{1/3}\sin \left(\dfrac{1}{3}\left(\dfrac{\pi}{2}+\theta \right)\right) + C_0, \\
\bar{\bfv}(x,y)   & = (10\alpha\sin^2 (\pi x) \sin(\pi y)\cos(\pi y), -10\alpha\sin^2 (\pi y) \sin(\pi x)\cos(\pi x) )^T,\\
\bar{q}(x,y)      & = \alpha \bar{p} (x,y),\\
\bar{\gamma}(x,y) & = \Pi_{[0,5]} \left (\dfrac{1}{\alpha}\bar{\bfu}\cdot\bar{\bfv} \right),
\end{align*}
where $(r,\theta)$ are the polar coordinates of $(x,y)$ and $C_0$ is a constant such that $\bar{p} \in L^2 _0 (\Omega)$.
Since $\alpha$ is very small, we use the same scheme as for $\PP_0$ case for the $\PP_1$ case. Although this scheme does not guarantee that $\gamma^h$ is bounded by $a$ and $b$, the numerical results show a good approximation. We summarize the number of iterations of each solver, the a priori error estimations with respect to the reference solutions and the a posteriori error estimations for our three discrete schemes in Figure \ref{error-ej2} and Tables \ref{ej2-p0}, \ref{ej2-p1} and \ref{ej2-semi}.  Due to $\bar{p}\in H^{4/3 - \varepsilon} (\Omega)$ for all $\varepsilon>0$ and a singularity of $\nabla\bar{p}$ at $(0,0)$, the theoretical convergence rates for the states, adjoints and control should not be reachable. However, we recover optimal rates obtained in Theorems \ref{taylor-hood}, \ref{conv-p0}, \ref{conv-p1} and \ref{semi} for the control and adjoints. In the case of the staes, we recover a convergence rate of $\mathcal{O}(h^{4/3})$. The effectivity index is asymptotic exact in $\PP_0$ and $\PP_1$ case and stable in every case.

\begin{table}[H]
\centering
\begin{tabular}
[c]{|c|c|c|c|c|c|c|c|}\hline
$h$ & \texttt{dof} & It. & $\Vert e(\gamma)\Vert_{0,\Omega}$ & $\Vert(e(\bfu),e(p))\Vert$ & $\Vert(e(\bfv),e(q))\Vert$ & $\eta$& $\theta$\\\hline
$1/4$   & $1128$   & $11$ & $6.85978E-01$ & $5.90663E-02$ & $6.45450E-04$ & $8.75376E-01$ & $1.271394$\\
$1/8$   & $4168$   & $12$ & $4.47430E-01$ & $2.56906E-02$ & $2.78855E-04$ & $4.87118E-01$ & $1.086912$\\
$1/16$  & $16008$  & $13$ & $2.40082E-01$ & $8.08743E-03$ & $7.25668E-05$ & $2.47780E-01$ & $1.031479$\\
$1/32$  & $62728$  & $14$ & $1.20819E-01$ & $2.30436E-03$ & $1.81712E-05$ & $1.21984E-01$ & $1.009459$\\
$1/64$  & $248328$ & $15$ & $6.05314E-02$ & $7.23438E-04$ & $4.54532E-06$ & $6.07436E-02$ & $1.003434$\\
$1/128$ & $988168$ & $16$ & $3.03072E-02$ & $2.50717E-04$ & $1.13645E-06$ & $3.03546E-02$ & $1.001530$\\\hline
\end{tabular}
\caption{A priori and a posteriori error estimates with effectivity indexes, $\PP_0$ case.}
\label{ej2-p0}
\end{table}

\begin{table}[H]
\centering
\begin{tabular}
[c]{|c|c|c|c|c|c|c|c|}\hline
$h$ & \texttt{dof} & It. & $\Vert e(\gamma)\Vert_{0,\Omega}$ & $\Vert(e(\bfu),e(p))\Vert$ & $\Vert(e(\bfv),e(q))\Vert$ & $\eta$& $\theta$\\\hline
$1/4$   & $1097$   & $11$ & $6.13171E-01$ & $1.02108E-01$ & $6.91900E-04$ & $7.90560E-01$ & $1.271784$\\
$1/8$   & $4009$   & $12$ & $2.05330E-01$ & $1.71637E-02$ & $2.32263E-04$ & $2.35752E-01$ & $1.144170$\\
$1/16$  & $15305$  & $13$ & $3.98752E-02$ & $4.10630E-03$ & $5.91178E-05$ & $4.84436E-02$ & $1.208488$\\
$1/32$  & $59785$  & $14$ & $1.23031E-02$ & $1.49056E-03$ & $1.48951E-05$ & $1.57407E-02$ & $1.270121$\\
$1/64$  & $236297$ & $15$ & $4.09242E-03$ & $5.73823E-04$ & $3.73297E-06$ & $5.55114E-03$ & $1.343303$\\
$1/128$ & $939529$ & $16$ & $1.45720E-03$ & $2.24958E-04$ & $9.34358E-07$ & $2.06812E-03$ & $1.402627$\\\hline
\end{tabular}
\caption{A priori and a posteriori error estimates with effectivity indexes, $\PP_1$ case.}
\label{ej2-p1}
\end{table}

\begin{table}[H]
\centering
\begin{tabular}
[c]{|c|c|c|c|c|c|c|c|}\hline
$h$ & \texttt{dof} & It. & $\Vert e(\gamma)\Vert_{0,\Omega}$ & $\Vert(e(\bfu),e(p))\Vert$ & $\Vert(e(\bfv),e(q))\Vert$ & $\eta$& $\theta$\\\hline
$1/4$   & $1032$   & $11$ & $4.62239E-02$ & $3.70425E-02$ & $5.21354E-04$ & $2.54485E-01$ & $4.296020$\\
$1/8$   & $3784$   & $12$ & $1.26431E-02$ & $1.14962E-02$ & $2.30288E-04$ & $7.78721E-02$ & $4.556619$\\
$1/16$  & $14472$  & $13$ & $2.13337E-03$ & $3.96809E-03$ & $5.91302E-05$ & $2.64559E-02$ & $5.871768$\\
$1/32$  & $56584$  & $14$ & $3.70431E-04$ & $1.48344E-03$ & $1.48965E-05$ & $9.74963E-03$ & $6.376201$\\
$1/64$  & $223752$ & $15$ & $8.18388E-05$ & $5.73347E-04$ & $3.73309E-06$ & $3.74673E-03$ & $6.469131$\\
$1/128$ & $889864$ & $16$ & $2.10000E-05$ & $2.25375E-04$ & $9.34795E-07$ & $1.47000E-03$ & $6.494269$\\
\hline
\end{tabular}
\caption{A priori and a posteriori error estimates with effectivity indexes, semi-discrete case.}
\label{ej2-semi}
\end{table}

\begin{figure}[H]
\centering
\includegraphics[height=4.2cm]{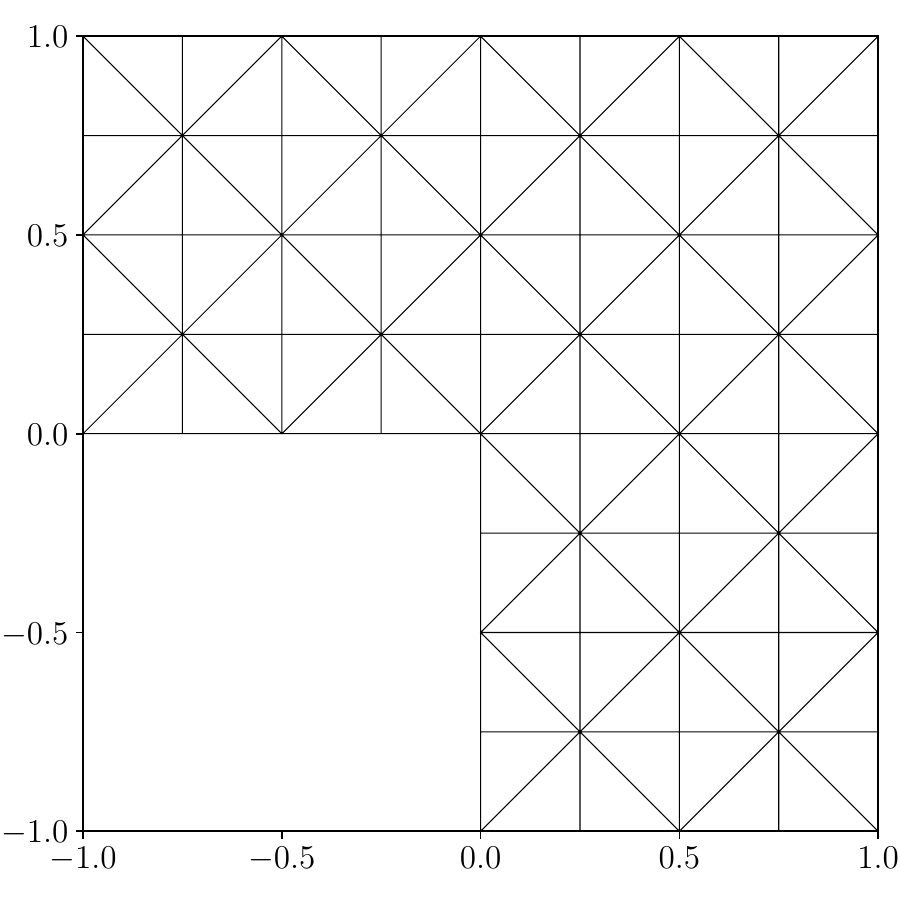}
\caption{Initial mesh for adaptive refinement}
\label{mesh-0}
\end{figure}

\begin{figure}[H]
\centering
\begin{subfigure}[b]{0.27\textwidth}
	 \centering
	 \includegraphics[width=\textwidth]{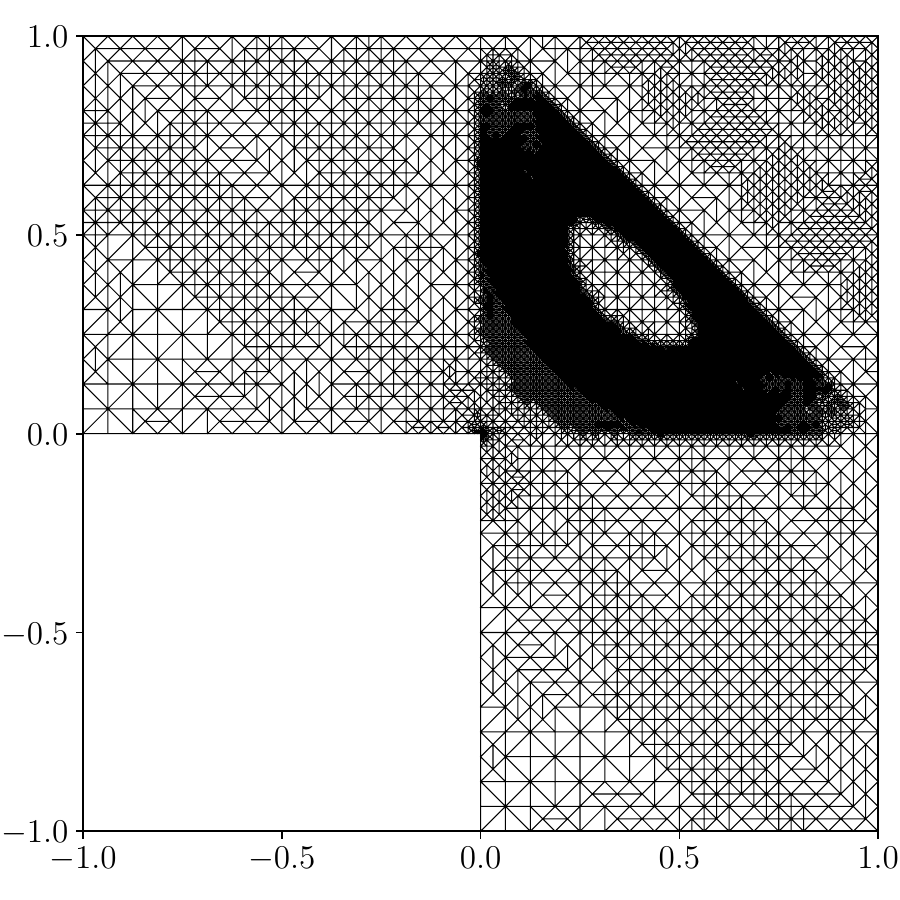}
\end{subfigure}
\qquad
\begin{subfigure}[b]{0.27\textwidth}
	 \centering
	 \includegraphics[width=\textwidth]{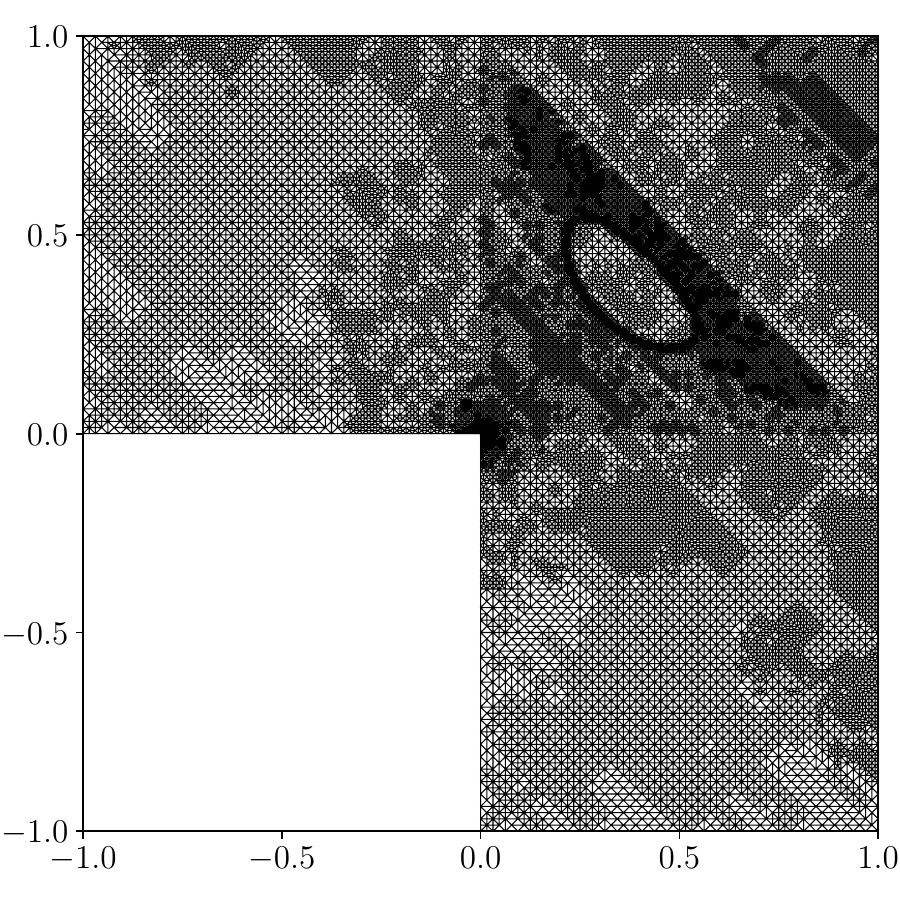}
\end{subfigure}
\qquad
\begin{subfigure}[b]{0.27\textwidth}
	 \centering
	 \includegraphics[width=\textwidth]{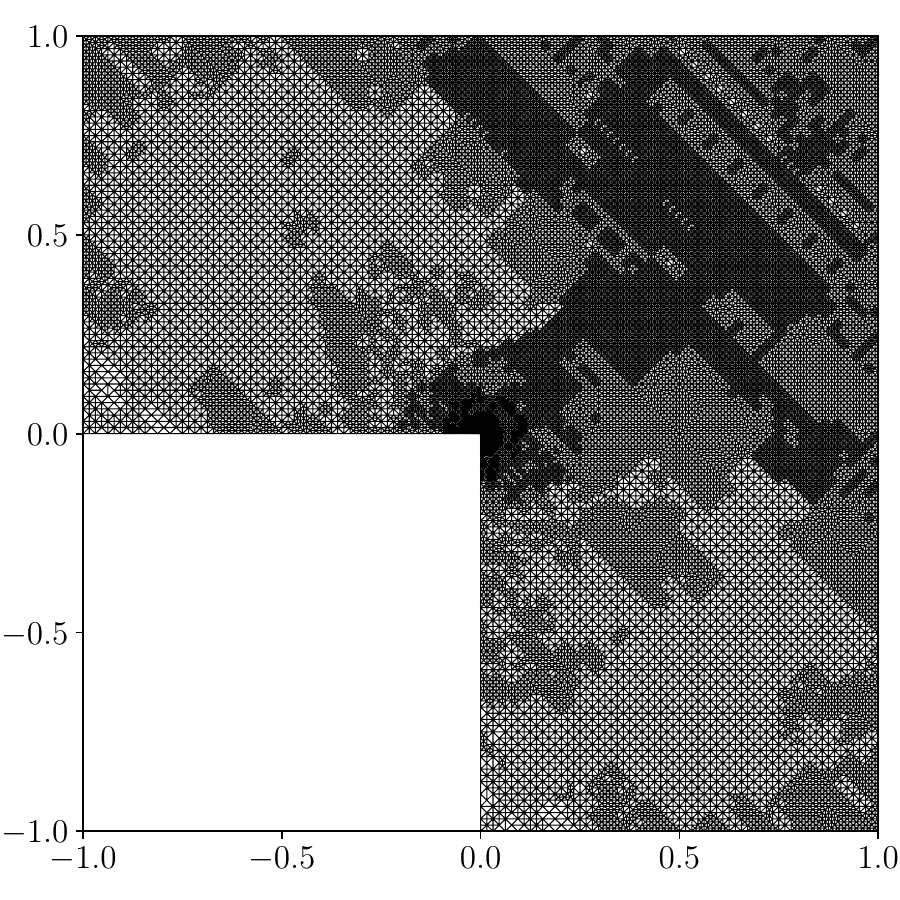}
\end{subfigure}
\caption{Adapted meshes for $\PP_0$ (left, stage 24, 497218 \texttt{dof}, 49634 elements), $\PP_1$ (center, stage 27, 440609 \texttt{dof}, 46119 elements) and semi-discrete schemes (right, stage 31, 451780 \texttt{dof}, 49898 elements).}
\label{meshes-ej2}
\end{figure}

\begin{figure}[H]
\centering
\begin{subfigure}[b]{0.3\textwidth}
	 \centering
	 \includegraphics[width=0.95\textwidth]{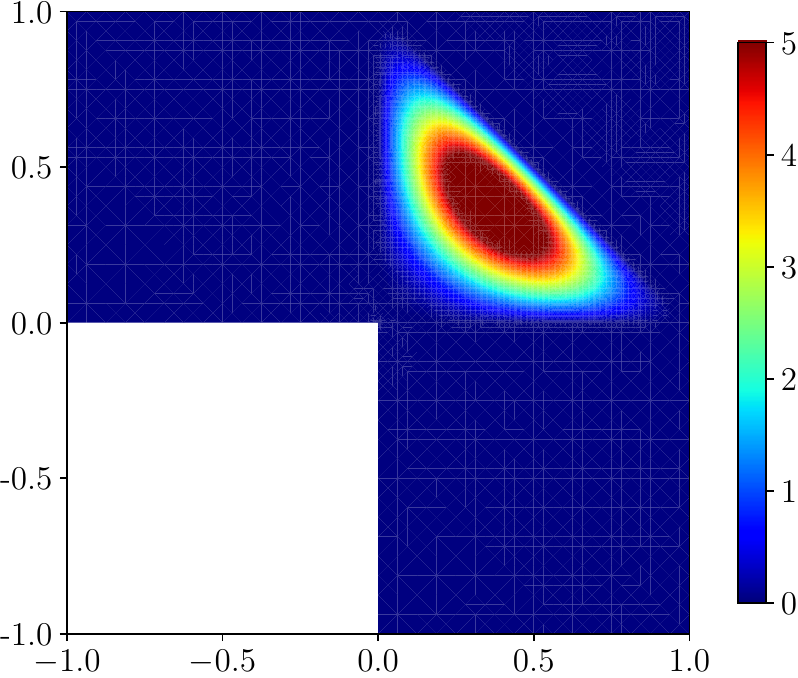}
\end{subfigure}
\quad
\begin{subfigure}[b]{0.3\textwidth}
	 \centering
	 \includegraphics[width=0.95\textwidth]{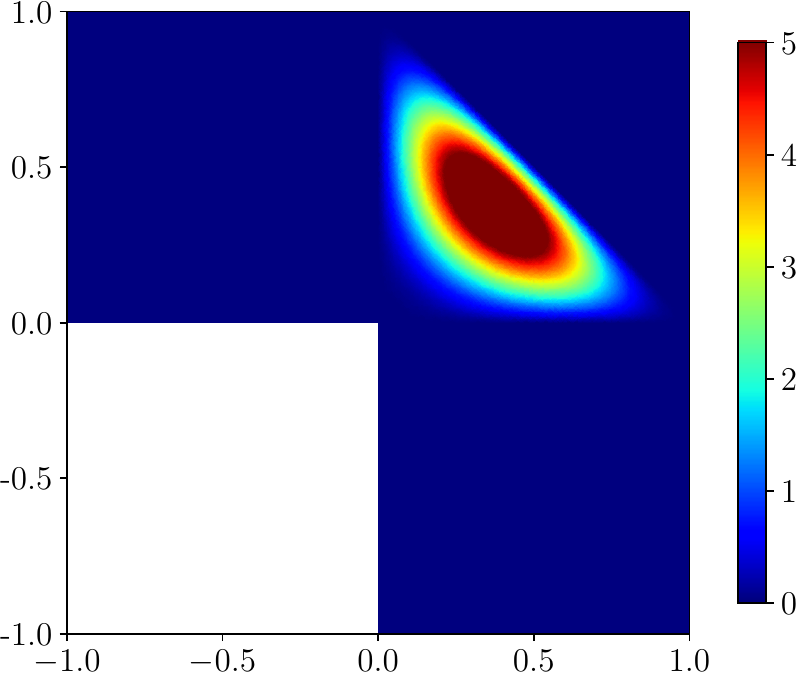}
\end{subfigure}
\quad
\begin{subfigure}[b]{0.3\textwidth}
	 \centering
	 \includegraphics[width=0.95\textwidth]{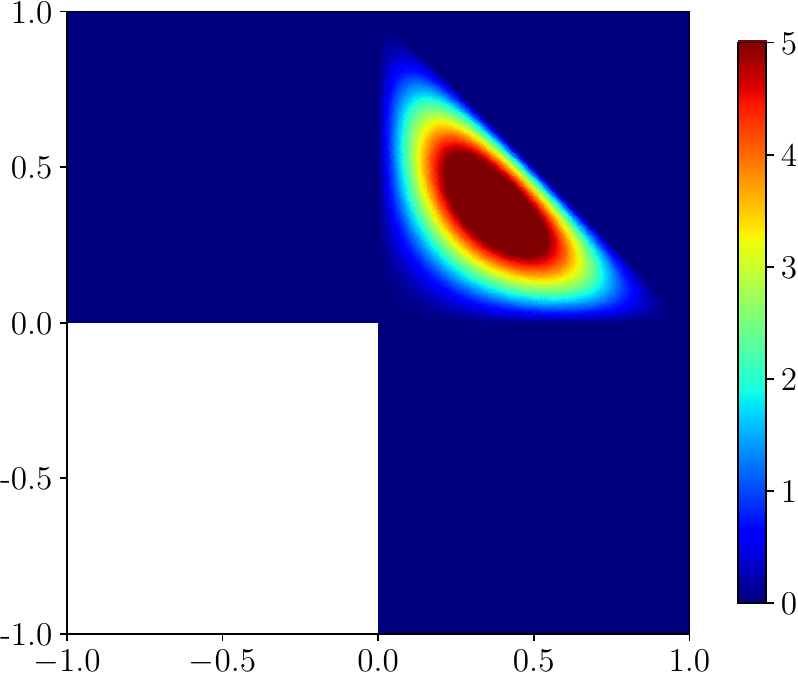}
\end{subfigure}
\caption{Isovalues of $\gamma^h$ on the adapted meshes for $\PP_0$, $\PP_1$  and semi-discrete schemes (from left to right).}
\label{gamma-ej2}
\end{figure}

\begin{figure}[H]
\centering
\includegraphics[width=0.9\textwidth]{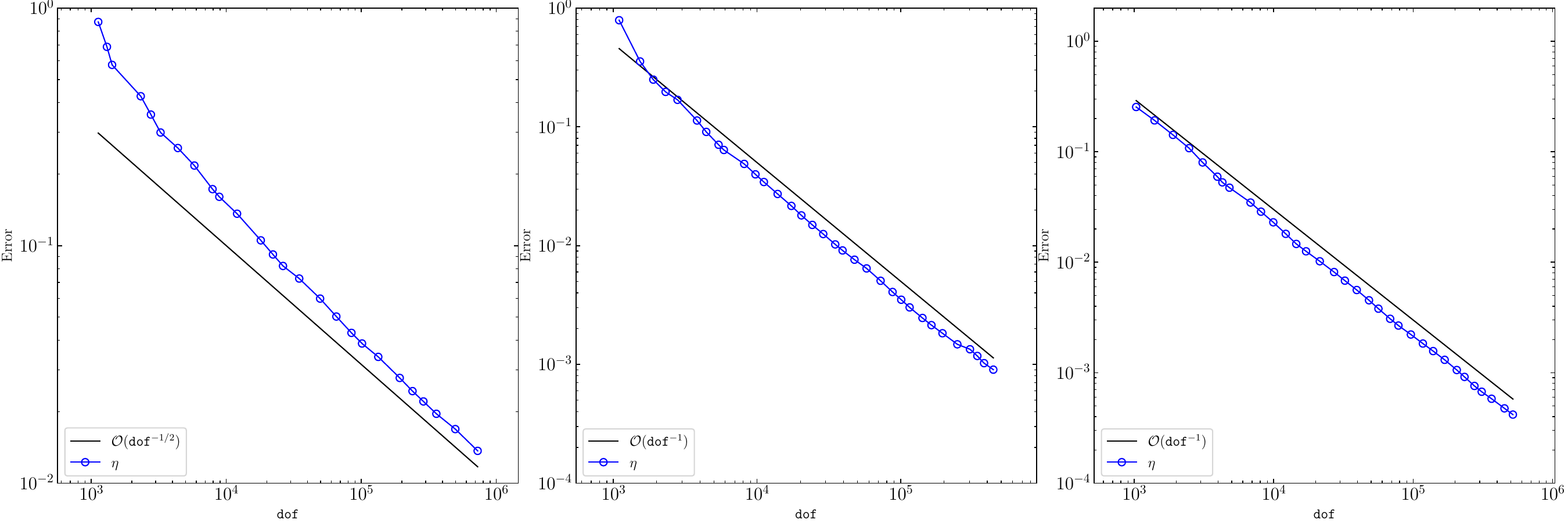}
\caption{History of a posteriori error estimator convergence for $\PP_0$ (left), $\PP_1$ (center) and semidiscrete schemes (right).}
\label{post-ej2}
\end{figure}

To improve the numerical solution, we perform the Algorithm \ref{afem} for $\rho=0.75$ and the initial mesh from Figure \ref{mesh-0}, and each discrete scheme, until the final refined mesh has more than $400000$ \texttt{dof}. The final meshes, plots of the optimal controls and the a posteriori error estimates are reported in Figures \ref{meshes-ej2}, \ref{gamma-ej2} and \ref{post-ej2}. In the discrete schemes, we observe that the adaptive refinement algorithm generates meshes with more elements near the point $(0,0)$ and the contour lines $\gamma=0$ and $\gamma=5$. However, this phenomenon is attenuated in the semi-discrete scheme since the a posteriori error estimate only considers the residual estimates for the forward and adjoint equations, recognizing the singularity of $\nabla\bar{p}$ at $(0,0)$. In all schemes, a posteriori estimator has a convergence order similar to the theoretical order.

\subsection{Adaptive refinement for obstacle recovering}
\begin{figure}[H]
\centering
\includegraphics[height=5cm]{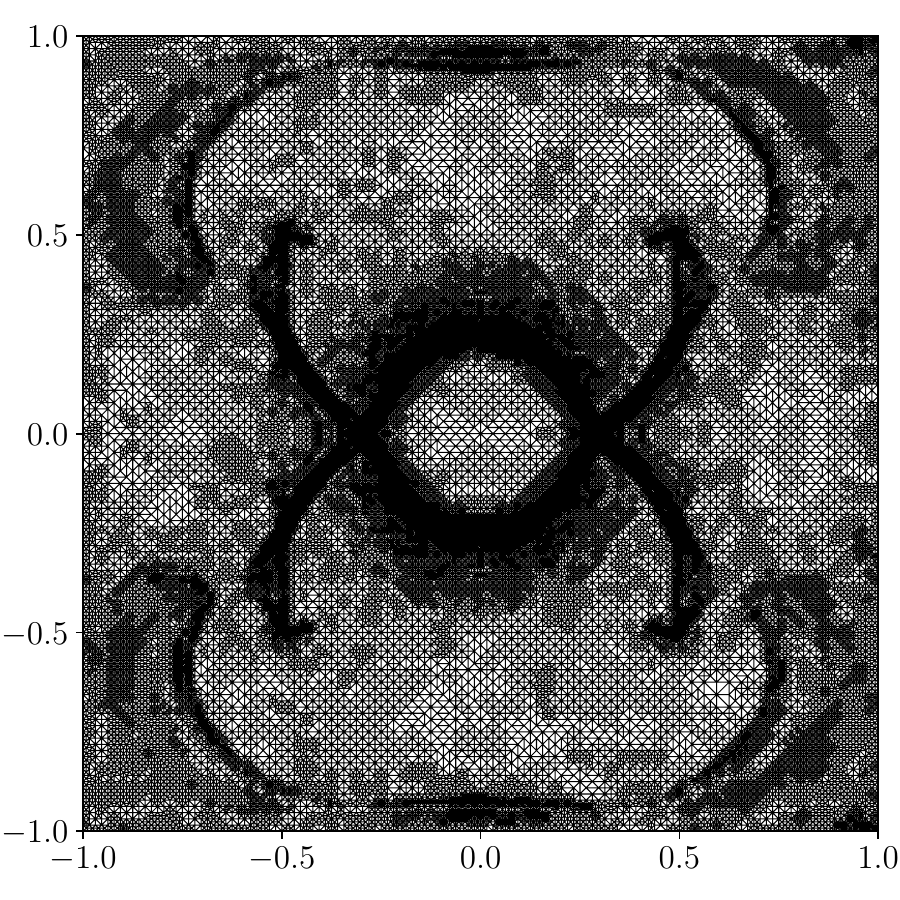}
\caption{Final mesh for adaptive refinement, stage 26, 107108 elements}
\label{mesh-obs}
\end{figure}

\begin{figure}[H]
\centering
\begin{subfigure}[c]{0.4\textwidth}
	 \centering
	 \includegraphics[height=5cm]{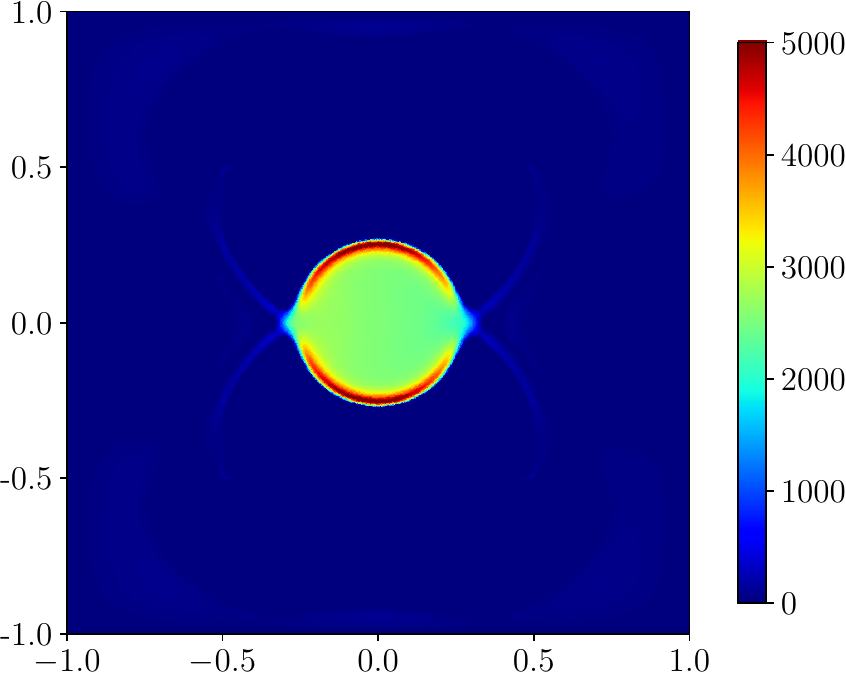}
\end{subfigure}
\qquad
\begin{subfigure}[c]{0.4\textwidth}
	 \centering
	 \includegraphics[height=5cm]{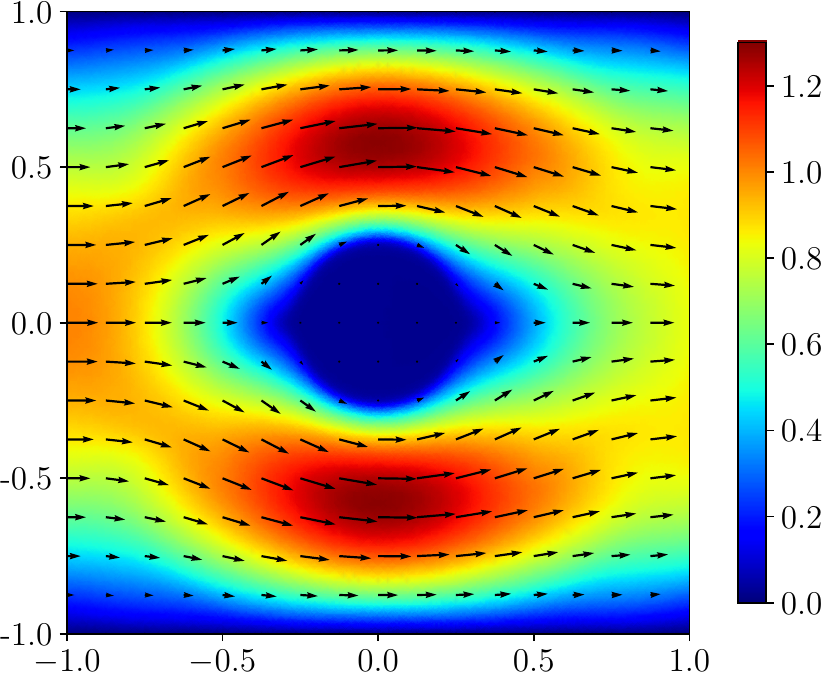}
\end{subfigure}
\caption{Plots of optimal $\gamma^h$ (left) and $\bfu^h$ (right), stage 26, 1056014 \texttt{dof}}
\label{opt-obs}
\end{figure}

\begin{figure}[H]
\centering
\includegraphics[height=5cm]{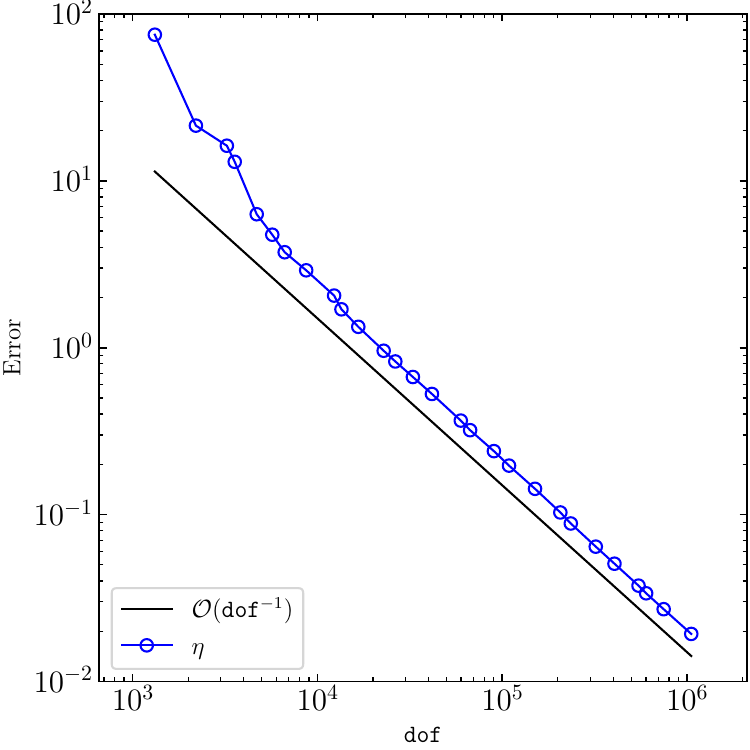}
\caption{History of the a posteriori error estimator convergence.}
\label{post-obs}
\end{figure}

In our final test, we seek to replicate the numerical results obtained in \cite[Section 4.2]{AO22} for the detection of an obstacle immersed in a fluid. First, in order to generate $\bfu_0$, we solve \eqref{VFh} in the domain $\Omega= ]-1,1[ ^2$ discretized by a hyperfine mesh (165440 elements), with the parameters $\bff=\boldsymbol{0}$, $\nu=1$, the boundary conditions
\begin{alignat*}{2}
\bfu  & =(1-y^{2},0)^{T} &\quad&\text{on }(\{-1\}\times [-1,1])\cup ([-1,1]\times\{-1,1\}),\\
(-\nu\nabla\bfu+p\boldsymbol{I})\bfn  & =\boldsymbol{0} && \text{on }\{1\}\times [-1,1],
\end{alignat*}
and $\gamma\in L^{2}(\Omega)$ given by
\[
\gamma=\left\{
\begin{array}
[c]{cl}%
10^{20} & \text{if }(x,y)\in B\\
0 & \text{otherwise}%
\end{array}
\right.
\]
where $B=\{(x,y)\in\Omega ~ \mid ~ (x^{2}+y^{2})^{1/2}\leq 0.25\}$. The set $B$ represents a circular obstacle immersed in the control volume $\Omega$.
\begin{figure}[H]
\centering
\includegraphics[height=5cm]{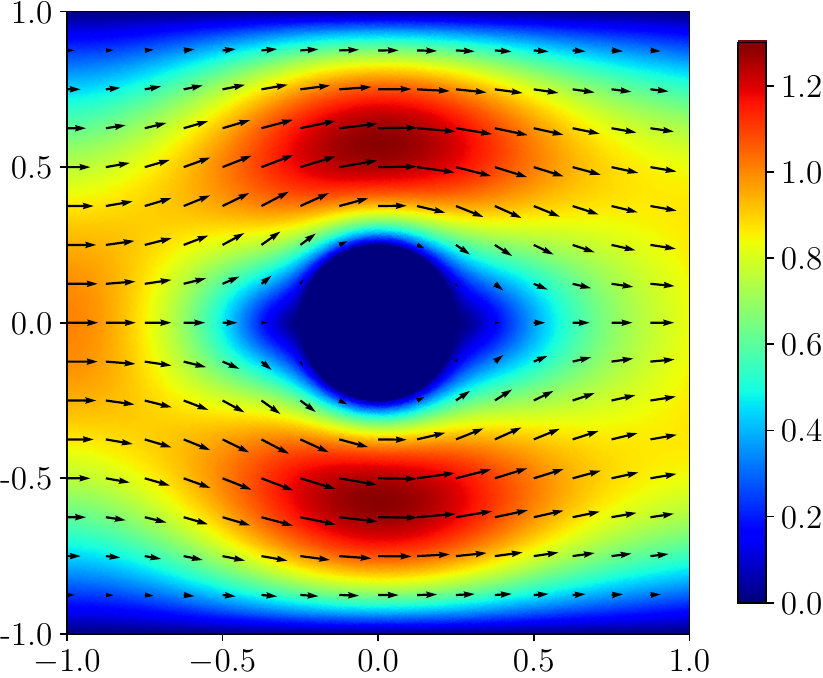}
\caption{Plot of $\bfu_0$ in a hyperfine mesh.}
\label{ref-obs}
\end{figure}
Then, we solve \eqref{min-h} with the semi-discrete scheme for $\alpha=10^{-10}$, $a=0$, $b=10^4$, $\bff=\boldsymbol{0}$, the measurement subset $\omega = ]-1/2,1/2[^2$, and Dirichlet boundary conditions given by $\bfu=\bfu_0$ on $\partial\Omega$. We perform the Algorithm \ref{afem} for $\rho=0.6$ and the initial mesh of Figure \ref{mesh-0} until the final refined mesh has more than $1000000$ \texttt{dof}. We report the final mesh, plots of the optimal control and optimal velocity, and the a posteriori error estimates in Figures \ref{mesh-obs}, \ref{opt-obs} and \ref{post-obs}. The adaptive refinement algorithm can refine the mesh near the boundary of the obstacle and the vertices of the measurement subset $\omega$, following the theoretical convergence order.

\section{Conclusion and future outlook}\label{S6}
In this paper, we have developed a comprehensive numerical framework to solve an optimal control problem governed by the steady Navier-Stokes-Brinkman equations, focusing on the identification of a scalar permeability parameter $\gamma$. For a discretization of 
the states and adjoint equations by Taylor-Hood finite elements, we deduced our a priori and a posteriori error estimates for two different fully discrete schemes for the control, using $\PP_0$ and $\PP_1$ finite elements, and for a semi-discrete scheme where the discrete control is given by an explicit expression that depends on the optimal states and adjoints.

Additionally, our reliable and efficient a posteriori error estimators are used to guide adaptive mesh refinement strategies, enabling the accurate recovery of smooth and discontinuous permeability parameters while significantly reducing computational costs. Numerical experiments validated the theoretical findings, demonstrating the effectiveness of the proposed method in practical scenarios, particularly with adaptive refinement improving accuracy near boundaries or discontinuities. The derivation of error estimates tailored to this problem and the development of adaptive refinement strategies for discontinuous parameters represent significant advancements. Furthermore, the semi-discrete scheme provides a flexible approach that bridges the gap between theoretical stability results and practical computational methods. 

Despite all the foundational contributions in this work, there is still room for improvement. An interesting post-processing method is discussed in \cite[Section 5.3]{KV09}, where an optimal control discretized in $\PP_1$ is projected to $\PP_0$ to improve the control convergence rate. Moreover, the performance of our solvers can be reduced when $\nu$ or $\alpha$ are small enough or $\gamma$ is large enough, which occurs in the resolution of realistic problems, resulting in the need for a stabilization method. In the case of equations \eqref{VFh} and \eqref{AVFh}, the SUPG \cite{TL91} and residual local projection \cite{ABPV12} schemes can contribute to stabilize the discrete optimal control problem. Furthermore, different penalty terms can be integrated into the cost function, such as using $L^1$ or the total variation norms \cite{DlR15}. In spite of the lack of differentiability of the cost function, the main advantage is to obtain a denoised state $\bfu$ from $\bfu_0$. By addressing these directions, future work can build upon the foundation laid in this paper, advancing the state-of-the-art in optimal control and parameter identification problems. 

\section*{Acknowledgments}
The authors acknowledge Axel Osses and Rodolfo Araya for the fruitful discussions about this article. Jorge Aguayo also thanks the funding of ANID CMM FB210005 Basal.

\bibliographystyle{elsarticle-harv}
\bibliography{ref.bib}

\end{document}